\documentclass{amsart}

\usepackage[a4paper,left=35mm,right=35mm,top=30mm,bottom=30mm,marginpar=25mm]{geometry}

\usepackage{amsfonts}
\usepackage{amsmath}
\usepackage{amssymb}
\usepackage{amsthm}
\usepackage[latin1]{inputenc}
\usepackage{eurosym}
\usepackage[dvips]{graphics}
\usepackage{graphicx}
\usepackage{epsfig}
\usepackage{hyperref}
\usepackage{dsfont}
\usepackage{color}
\usepackage{xcolor}
\usepackage{booktabs, multirow} % for borders and merged ranges
\usepackage{soul}% for underlines
\usepackage{changepage,threeparttable} % for wide tables
\usepackage{graphicx,subcaption}
\usepackage[displaymath,mathlines]{lineno}

\allowdisplaybreaks

\usepackage{hyperref}
\usepackage{mathtools}
\usepackage{ifthen}
\graphicspath{./figures/}
\makeindex

\newcommand{\bK}{\mathbf{K}}

\newcommand{\Om}{\Omega}
\newcommand{\R}{\mathbb{R}}
\newcommand{\Ll}{\mathcal{L}}
\newcommand{\rest}{%
	\,\raisebox{-.127ex}{\reflectbox{\rotatebox[origin=br]{-90}{$\lnot$}}}\,%
}

\numberwithin{equation}{section}

\newtheoremstyle{thmlemcorr}{10pt}{10pt}{\itshape}{}{\bfseries}{.}{10pt}{{\thmname{#1}\thmnumber{
			#2}\thmnote{ (#3)}}}
\newtheoremstyle{thmlemcorr*}{10pt}{10pt}{\itshape}{}{\bfseries}{.}\newline{{\thmname{#1}\thmnumber{
			#2}\thmnote{ (#3)}}}
\newtheoremstyle{defi}{10pt}{10pt}{\itshape}{}{\bfseries}{.}{10pt}{{\thmname{#1}\thmnumber{
			#2}\thmnote{ (#3)}}}
\newtheoremstyle{remexample}{10pt}{10pt}{}{}{\bfseries}{.}{10pt}{{\thmname{#1}\thmnumber{
			#2}\thmnote{ (#3)}}}
\newtheoremstyle{ass}{10pt}{10pt}{}{}{\bfseries}{.}{10pt}{{\thmname{#1}\thmnumber{
			A#2}\thmnote{ (#3)}}}

\theoremstyle{thmlemcorr}
\newtheorem{theorem}{Theorem}
\numberwithin{theorem}{section}

\theoremstyle{thmlemcorr*}
\newtheorem{theorem*}{Theorem}
\newtheorem{lemma*}[theorem]{Lemma}
\newtheorem{corollary*}[theorem]{Corollary}
\newtheorem{proposition*}[theorem]{Proposition}
\newtheorem{problem*}[theorem]{Problem}
\newtheorem{conjecture*}[theorem]{Conjecture}

\theoremstyle{defi}
\newtheorem{definition}[theorem]{Definition}

\theoremstyle{remexample}
\newtheorem{remark}[theorem]{Remark}

\theoremstyle{ass}

\begin{document}
	
	\title[Computational methods for nonlocal mean field games with applications]{Computational methods for nonlocal mean field games with applications}
	
	%\author{Osher Lab}
	
	\author{Siting Liu}
	\address[S. L.]{Department of Mathematics, UCLA}
	\email{siting6@math.ucla.edu}
	
	\author{Matthew Jacobs}
	\address[M. J.]{Department of Mathematics, UCLA}
	\email{majaco@math.ucla.edu}
	
	\author{Wuchen Li}
	\address[W. L.]{Department of Mathematics, UCLA}
	\email{wcli@math.ucla.edu}
	
	\author{Levon Nurbekyan}
	\address[L. N.]{Department of Mathematics, UCLA}
	\email{lnurbek@math.ucla.edu}
	
	\author{Stanley J. Osher}
	\address[S. J. O.]{Department of Mathematics, UCLA}
	\email{sjo@math.ucla.edu}
	
	
	%These are the 2020 classifications.
	%35A15 Variational methods applied to PDEs
	%35Q89 PDEs in connection with mean field game theory
	%35Q91 PDEs in connection with game theory, economics, social and behavioral sciences
	%35Q93 PDEs in connection with control and optimization
	%49N80 Mean field games and control
	%91A16 Mean field games (aspects of game theory)
	%65M70 Spectral, collocation and related methods for ini- tial value and initial-boundary value problems in- volving PDEs
	%93A15 Large-sacels systems
	%93A16 Multi-agent systems

	\thanks{This work was supported by AFOSR MURI FA9550-18-1-0502, AFOSR Grant No. FA9550-18-1-0167, ONR Grant No. N00014-18-1-2527, ONR N00014-20-1-2093. L.N. was also supported by Simons Foundation and the Centre de Recherches Math\'{e}matiques, through the Simons-CRM scholar-in-residence program.}
	\date{\today}
	
	\begin{abstract}
		We introduce a novel framework to model and solve mean-field game systems with nonlocal interactions extending the results in \cite{nursaude18}. Our approach relies on kernel-based representations of mean-field interactions and feature-space expansions in the spirit of kernel methods in machine learning. We demonstrate the flexibility of our approach by modeling various interaction scenarios between agents. Additionally, our method yields a computationally efficient saddle-point reformulation of the original problem that is amenable to state-of-the-art convex optimization methods such as the primal-dual hybrid gradient method (PDHG). We also discuss potential applications of our methods to multi-agent trajectory planning problems.
	\end{abstract}
	
	\maketitle
	
	%\tableofcontents

	\section{Introduction}
	
	In this paper, we study computational and modeling aspects of nonlocal mean-field game systems. Specifically, we are interested in the system
	\begin{equation} \label{eq:main_gen}
	\begin{cases}
	-\phi_t  + H(x,\nabla \phi) = f\left(x, \int_{\Om} K(x,y) \rho(y,t) dy \right) \\
	\rho_t  -\nabla \cdot (\rho \nabla_p H(x,\nabla \phi))=0\\
	\rho(x,0)=\rho_0(x),~\phi(x,1)=g\left(x,\int_{\Om} S(x,y) \rho(y,1) dy\right)
	\end{cases}
	\end{equation}
	Above, $\Omega\subset \R^d$ is a domain with smooth boundary, $H:\Om \times \R^d \to \R$ is the Hamiltonian, $K,S:\Om\times \Om \to \R$ are interaction kernels, and $f,g:\Omega\times \R \to \R$ are interaction functions.
	
	System \eqref{eq:main_gen} describes Nash equilibria in an infinite-dimensional differential game where a continuum of agents interact through their distribution in the state-space. These games are called mean-field games (MFG) and were introduced by Huang, Malham\'{e}, and Caines \cite{HCM06,HCM07}, and Lasry and Lions \cite{LasryLions06a,LasryLions06b,LasryLions2007}. Currently, MFG is a thriving research direction with applications in economics \cite{gueantlasrylions11,moll14,Gomes:2015th,moll17}, finance \cite{lehalle16,lehalle18,jaimungal19,caines17}, crowd motion \cite{wolfram11,wolfram14,aurell18,achdou19}, industrial engineering \cite{kizikale19,paola19,gomessaude'18},  data science \cite{EHanLi2018_meanfield}, material dynamics \cite{welch19}, and more. For a detailed introduction and review of MFG theory we refer to \cite{LasryLions2007,gueantlasrylions11,CardaNotes,gomsaude'14}.
	
	In \eqref{eq:main_gen}, $\rho(\cdot,t)$ represents the population density at time $t$, and $\rho_0$ is an initial distribution. Furthermore, $(x,t)\mapsto \phi(x,t)$ is a value function that measures the optimal value of an agent at position $x$ and time $t$. Functions $f,g$ and kernels $K,S$ model the interaction between a single agent and the population. 
	
	Throughout the paper, we assume that $f,g,K,S$ are $C^2$ functions, and a pair $(\phi,\rho)$ is a weak solution of \eqref{eq:main_gen} if $\phi$ is a viscosity and $\rho$ is a distribution solution for HJB and continuity equations, respectively. The PDE theory of \eqref{eq:main_gen} is well understood in this setting, and we refer to \cite{LasryLions2007,CardaNotes} for a detailed exposition of the subject.
	
	Our goal is to develop computational and modeling methods for \eqref{eq:main_gen}. There are several general purpose numerical methods for MFG that apply to \eqref{eq:main_gen}. In \cite{achdou10,achdou12,achdou13,achdou13b}, the authors develop and analyze finite-difference methods for \eqref{eq:main_gen} and related models. Their approach is based on a solution of the discrete problem applying Newton's method.
	
	For so-called \textit{potential} MFG systems, there are several primal-dual convex optimization methods such as alternating direction method of multipliers (ADMM) \cite{bencar'15,bencarsan'17} and primal-dual hybrid gradient (PDHG) \cite{silva18,silva19}. These methods are extensions of the celebrated Brenier-Benamou method for computation of optimal transport maps \cite{BenamouBrenier2000} to the MFG setting.
	
	In this paper, our goal is to develop computational and modeling methods specifically for nonlocal MFG systems. The term nonlocal refers to expressions
	\begin{equation*}
	f\left(x, \int_{\Om} K(x,y) \rho(y,t) dy \right),\quad g\left(x,\int_{\Om} S(x,y) \rho(y,1) dy\right).
	\end{equation*}
	Indeed, the calculation of these terms at $x$ requires the knowledge of $\rho(y,t)$ at all $y$. For general MFG systems, these terms are replaced by $f(x,\rho),~g(x,\rho)$ where one allows local interactions as well; that is, $f(x,\rho)=\tilde{f}(x,\rho(x,t))$, and $g(x,\rho)=\tilde{g}(x,\rho(x,1))$ where, $\tilde{f}, \tilde{g}:\Omega\times \R_{+}\to \R$ are some functions, and $\rho(x,t)$ is the density of $\rho$. Note that for local interactions, the calculation of $f,~g$ terms at $x$ requires information only at $x$.
	
	The majority of numerical methods above apply to general models, including nonlocal interactions. However, they suit best the local ones, where a discretization of $\rho$ on a grid yields a straightforward discretization for interactions. For the nonlocal case though, the calculation of interaction terms requires matrix multiplication on a full grid to evaluate the expressions $\int_\Omega K(x,y) \rho(y,t)dy,~\int_\Omega S(x,y) \rho(y,1)dy$. Our approach solves this problem by encoding the interactions in a small number of \textit{expansion coefficients}.
	
	Furthermore, a critical feature of primal-dual methods mentioned above is that one of the proximal steps results in a decoupled system of one-dimensional convex optimization problems at the grid-points. Therefore, this step is parallelizable and yields a linear computational cost in the number of grid-points. However, direct applications of these methods to nonlocal problems yield fully coupled systems that are not parallelizable and yield a superlinear computational cost. Our method solves this problem as well. The expansion coefficients, that encode the interactions, decouple the aforementioned systems. Furthermore, we update these coefficients by an explicit proximal step that yields a linear computational cost.
	
	Our method is also well-suited for the Lagrangian framework. Indeed, in \cite{nursaude18}, where this approach was introduced, the authors solved \eqref{eq:main_gen} in Lagrangian coordinates. Thus, this approach paves a way to efficient computational methods for high-dimensional MFG problems.
	
	Another appealing feature of our method is the flexibility of modeling interactions. We expand $K,S$ in a basis that can be interpreted as features from kernel methods in ML \cite[Chapter 6]{mohri}. This allows us to design various interactions by only manipulating the basis. In particular, we can easily model heterogeneous regimes where agents interact only within specific subdomains and other interesting scenarios. Additionally, we can handle nonlocal interactions that are given by differential operators \cite[Tests 5, 6]{achdou10}.
	
	Finally, we would like to point out potential applications of our methods to multi-agent trajectory-planning. In general, even single-agent trajectory-panning problems are highly complex. With the number of agents increasing in a system, problems quickly become computationally overwhelming.
	
	A critical difficulty comes from modeling and computing the interactions. MFG theory provides a flexible framework to solve this problem. Theoretically, MFG solutions are optimal only when the number of agents is infinite. Nevertheless, one can generate sub-optimal trajectories that have appealing properties such as no-collision. Since our method provides a way of encoding mean-field interactions in a few coefficients independent of the number of agents, it is potentially applicable to large multi-agent trajectory planning problems.
	
	The paper has the following organization. In Section \ref{sec:method}, we provide a detailed description of our method. In Section \ref{sec:pdhg}, we devise a PDHG algorithm to solve \eqref{eq:main_gen} based on our method and discuss the computational efficiency of our approach. In Section \ref{sec:kernels}, we show how to model and approximate interactions using kernel methods from ML. Next, in Section \ref{sec:multi_agent_control}, we discuss potential applications of our methods to multi-agent trajectory planning problems. Finally, in Section \ref{sec:numerical_examples} we provide several numerical experiments.
	
	\section{The method of coefficients}\label{sec:method}
	
	For the sake of simplicity, we assume that $f(x,\zeta)=\zeta$, and $g(x,\zeta)=g(x)$; that is, we consider the system
	\begin{equation} \label{eq:main}
	\begin{cases}
	-\phi_t  + H(x,\nabla \phi) = \int_{\Om} K(x,y) \rho(y,t) dy \\
	\rho_t  -\nabla \cdot (\rho \nabla_p H(x,\nabla \phi))=0\\
	\rho(x,0)=\rho_0(x),~\phi(x,1)=g(x)
	\end{cases}
	\end{equation}
	The essence of the method is in an expansion of $K(x,y)$ in a family of functions. More precisely, assume that $\{f_i\}_{i=1}^r \subset C^2(\Omega)$ is an arbitrary family of functions. Furthermore, suppose that
	\begin{equation}\label{eq:K_expansion}
	K(x,y)=\sum_{i,j=1}^r k_{ij}f_i(x) f_j(y)
	\end{equation}
	\begin{remark}
	In general, $K$ may not have the form \eqref{eq:K_expansion}. In such cases, we approximate $K$ with kernels of such form.
	\end{remark}
	We denote by $\bK=(k_{ij})\in M_{r\times r}(\R)$, and, without loss of generality, we assume that $\bK$ is invertible. A straightforward calculation yields that
	\begin{equation*}
	\begin{split}
	\int_\Omega K(x,y)\rho(y,t)dt=\sum_{i=1}^r a_i(t) f_i(x),
	\end{split}
	\end{equation*}
	where
	\begin{equation*}
	a_i(t)=\sum_{j=1}^r k_{ij} \int_\Omega f_j(y)\rho(y,t)dy.
	\end{equation*}
	This means that no matter what $\rho$ is, the expression $\int_\Omega K(x,y)\rho(y,t)dy$ is a combination of $\{f_i\}_{i=1}^r$ with some \textit{unknown coefficients} $a_1,a_2,\cdots,a_r$. These coefficients encode all necessary information about the interactions. Thus, $(a_i)$ will be our new unknowns. Note that once we find $(a_i)$ we can solve \eqref{eq:main} by solving decoupled HJB and continuity equations.
	
	It turns out that $(a_i)$ are zeroes of a specific operator that is \textit{monotone} if $\rho \mapsto \int_\Omega K(x,y)\rho(y)dy$ is monotone, and a \textit{gradient} when $K$ is symmetric. To state the results precisely, denote by $\phi_a$ the viscosity solution of the HJB
	\begin{equation}\label{eq:HJB_a}
	\begin{cases}
	-\phi_t  + H(x,\nabla \phi) = \sum_{i=1}^r a_i(t) f_i(x) \\
	\phi(x,1)=g(x)
	\end{cases}
	\end{equation}
	Furthermore, denote by $\rho_a$ the distributional solution of the continuity equation
	\begin{equation}\label{eq:ce_a}
	\begin{cases}
	\rho_t  -\nabla \cdot(\rho \nabla_p H(x,\nabla \phi_a))=0\\
	\rho(x,0)=\rho_0(x)
	\end{cases}
	\end{equation}
	The following two theorems are the basis of our approach.
	\begin{theorem}{\cite[Theorem 2.3]{nursaude18}}\label{thm:convexity}
		The functional $a\mapsto \int_\Omega \phi_a(x,0)\rho_0(x)dx$ is concave and everywhere Fr\'{e}chet differentiable. Moreover,
		\begin{equation}\label{eq:dGa}
		\frac{\delta}{\delta a_i}\int_\Omega \phi_a(x,0)\rho_0(x)dx =\int_\Omega f_i(x) \rho_a(x,\cdot) dx,~1\leq i \leq r.
		\end{equation}
	\end{theorem}
	\begin{theorem}{\cite[Theorem 3.1]{nursaude18}}\label{thm:optimization}
		\begin{itemize}
			\item[i.] A pair $(\phi,\rho)$ is a solution of \eqref{eq:main} if and only if $(\phi,\rho)=(\phi_{a},\rho_{a})$ for some $a\in C\left([0,1];\mathbb{R}^r\right)$ such that
			\begin{equation}\label{eq:a_fixedpoint}
			a=\mathbf{K} \frac{\delta}{\delta a} \int_\Omega \phi_a(x,0)\rho_0(x)dx
			\end{equation}
			
			\item[ii.] If $\mathbf{K}$ is monotone, then \eqref{eq:a_fixedpoint} is equivalent to finding a zero of a monotone operator $a\mapsto \mathbf{K}^{-1} a - \frac{\delta}{\delta a} \int_\Omega \phi_a(x,0)\rho_0(x)dx,~a\in C\left([0,1];\mathbb{R}^r\right)$.
			
			\item[iii.] Additionally, if $\mathbf{K}$ is symmetric, \eqref{eq:a_fixedpoint} is equivalent to the convex optimization problem
			\begin{equation}\label{eq:a_equation_optimization}
			\begin{split}
			&\inf_{a\in C\left([0,1];\mathbb{R}^r\right)}  \frac{\langle \mathbf{K}^{-1} a, a \rangle }{2}- \int_\Omega \phi_a(x,0)\rho_0(x)dx,
			\end{split}
			\end{equation}
			where $\langle a,b\rangle=\sum_{i=1}^r \int_0^1  a_i(t)b_i(t)dt$ for $a,b\in C\left([0,1];\mathbb{R}^r\right)$.
		\end{itemize}
	\end{theorem}
	\begin{remark}
		Several remarks are in order.
		\begin{enumerate}
			\item Theorem \ref{thm:optimization} asserts that instead of finding $(\phi,\rho)$ in \eqref{eq:main} we just need to find the right coefficients $(a_i)$ and then solve the decoupled equations \eqref{eq:HJB_a}, \eqref{eq:ce_a}. 
			\item Here, we do not concentrate on technical aspects of Theorems \ref{thm:convexity}, \ref{thm:optimization} that are available for a periodic case, $\Omega=\mathbb{T}^d$, in \cite{nursaude18}. For the non-periodic case, \eqref{eq:main} and \eqref{eq:HJB_a} must be complemented with the no-flux Neumann boundary condition $\nabla_p H(x,\nabla \phi) \cdot \nu=0$ on $\partial \Omega$.
			\item $\bf K$ is monotone if and only if $\rho \mapsto \int_\Omega K(x,y)\rho(y)dy$ is monotone; that is, $\int_{\Omega^2} K(x,y)(\rho_2(x)-\rho_1(x))(\rho_2(y)-\rho_1(y))dxdy \geq 0$ for all $\rho_1,\rho_2\in \mathcal{P}(\Omega)$. This condition is essential for the uniqueness of solutions of \eqref{eq:main} \cite{LasryLions2007}.
			\item $\bf K$ is symmetric if and only if $K$ is symmetric; that is, $K(x,y)=K(y,x)$ for all $x,y\in \Omega$.
			\item For monotone interactions, $(a_i)$ are solutions of a monotone variational inequality or a convex optimization problem. Therefore, one can apply powerful convex optimization techniques to find $(a_i)$.
			\item Equation \eqref{eq:dGa} is critical for update rules of $(a_i)$. Indeed, it gives the ascent direction $\int_\Omega \phi_a(x,0)\rho_0(x)dx$ with respect to $(a_i)$. Moreover, this direction depends on $\rho_a$, which is available once $\phi_a$ is (approximately) computed at current $(a_i)$. As we shall see below, this property yields extremely simple update rules for $(a_i)$.
			\item These previous points are also very appealing for potential applications of the methods here to multi-agent trajectory planning problems.
			\item Optimization problem \eqref{eq:a_equation_optimization} is equivalent to the infinite-dimensional optimal control problem (58)-(59) in \cite{LasryLions2007}.
			\item There is nothing special about the choice $f(x,\zeta)=\zeta,~g(x,\zeta)=g(x)$ except simplicity. Analogous results can be proven for general $f,g$ as well. This topic is a subject of our future research.
			\item Precisely the same theorems are valid for stochastic systems as well. Again, we will address this topic in our future work.
		\end{enumerate}
	\end{remark}
	
\section{A primal-dual hybrid gradient algorithm}\label{sec:pdhg}

Here, we formulate \eqref{eq:a_equation_optimization} as a convex-concave saddle point problem and devise a primal-dual hybrid gradient (PDHG) algorithm of Chambolle and Pock \cite{champock11,champock16} to solve it. Introducing a Lagrange multiplier for the HJB equation in \eqref{eq:HJB_a} we obtain
\begin{equation}\label{eq:infsup_derivation}
\begin{split}
&\inf_{a\in C\left([0,1];\mathbb{R}^r\right)}  \frac{\langle \mathbf{K}^{-1} a, a \rangle }{2}- \int_\Omega \phi_a(x,0)\rho_0(x)dx\\
=&\inf_{\phi,a} \left\{ \frac{\langle \mathbf{K}^{-1} a, a \rangle }{2}  -\int_{\Omega}\phi(x,0)\rho_0(x) dx~\mbox{s.t.} ~ \eqref{eq:HJB_a} ~ \mbox{holds} \right\}\\
= & \inf_{\substack{\phi(x,1) =g\\ a}} \sup_{\rho}   \Bigg\{ \frac{\langle \mathbf{K}^{-1} a, a \rangle }{2} -\int_{\Omega}\phi(x,0)\rho_0(x) dx \\
& + \int_{\Omega}\int_0^1 \rho\left(	-\phi_t  + H(x,\nabla \phi) - \sum_{i=1}^r a_i(t) f_i(x) \right) dt dx \Bigg\}\\
=& \inf_{\substack{\phi(x,1) =g\\ a}} \sup_{\rho,v}   \Bigg\{ \frac{\langle \mathbf{K}^{-1} a, a \rangle }{2} -\int_{\Omega}\phi(x,0)\rho_0(x) dx - \int_{\Omega}\int_0^1 \left(\rho \phi_t  + \rho v\cdot \nabla \phi \right) dt dx\\
 &-\int_{\Omega}\int_0^1 \rho\left( L(x,v) + \sum_{i=1}^r a_i(t) f_i(x) \right) dt dx \Bigg\},\\
\end{split}
\end{equation} 
where we used the convex duality
\begin{equation*}
H(x,\nabla \phi)=\sup_v -v\cdot \nabla \phi - L(x,v).
\end{equation*}
The saddle point problem in \eqref{eq:infsup_derivation} is convex in $(\phi,a)$ and almost concave in $(\rho,v)$. Non-concavity comes from the terms $\rho v \cdot \phi$ and $\rho L(x,v)$. Following \cite{BenamouBrenier2000}, we remedy this problem by replacing $v$ with a flux variable $m=\rho v$. Thus, we obtain an equivalent saddle point problem
\begin{equation}\label{eq:infsupL}
\begin{split}
& \inf_{\substack{\phi(x,1) =g\\ a}} \sup_{\rho,m}   \Bigg\{ \frac{\langle \mathbf{K}^{-1} a, a \rangle }{2} -\int_{\Omega}\phi(x,0)\rho_0(x) dx - \int_{\Omega}\int_0^1 \left(\rho \phi_t  + m\cdot \nabla \phi \right) dx dt\\
&-\int_{\Omega}\int_0^1 \rho\left( L\left(x,\frac{m}{\rho}\right) + \sum_{i=1}^r a_i(t) f_i(x) \right) dx dt \Bigg\}\\
=&\inf_{\substack{\phi(x,1) =g\\ a}} \sup_{\rho,m}~\mathcal{L}(\phi,a,\rho,m)
\end{split}
\end{equation}
Note that $\left(\phi,a\right) \mapsto \mathcal{L}(\phi,a,\rho,m)$ is convex,  $\left(\rho,m\right) \mapsto \mathcal{L}(\phi,a,\rho,m)$ is concave, and the coupling between $\left(\phi,a\right)$ and $\left(\rho,m \right)$ is bilinear. Thus, we can apply PDHG \cite{champock11,champock16} to solve \eqref{eq:infsupL}. Also, note that the first-order optimality conditions for $\phi$ are
\begin{equation*}
\begin{cases}
\rho_t+\nabla\cdot m =0 \quad \mbox{in} \quad \Omega\times(0,1)\\
m(x,t) \cdot \nu=0 \quad \mbox{in} \quad  \partial \Omega\times (0,1)\\
\rho(x,0)=\rho_0(x) \quad \mbox{in} \quad  \Omega
\end{cases}
\end{equation*}
Therefore, the no-flux boundary condition for $m$ and the initial condition for $\rho$ are incorporated in \eqref{eq:infsupL}, and no extra considerations are necessary.

Furthermore, we must add a constraint $\rho\geq 0$ to account for $\rho$ being a probability distribution. This adjustment is coherent with the derivation \eqref{eq:infsup_derivation} because the viscosity-solution constraint \eqref{eq:HJB_a} should be replaced by pointwise constraints $-\partial_t+H(x,\nabla \phi)\leq \sum_{i=1}^r a_i(t) f_i(x),~\phi(x,1)\leq g(x)$. We refer to \cite{bencarsan'17} for details.

As illustrated in \cite{matt19,mattflavien}, the choices of spaces for variables are crucial when applying PDHG. Correct choices render algorithms with grid-size-independent convergence rates. For $a,\rho,m$ we choose $L^2$ spaces, whereas for $\phi$ we choose $H^1$. The motivation for this choice comes from convergence analysis of PDHG. Indeed, upper bounds on the step-sizes depend on the inverse of the norm of the bilinear coupling
\begin{equation*}
	\left|\int_{\Omega}\int_0^1 \left(\rho \phi_t  + m\cdot \nabla \phi \right) dx dt  \right| \leq \|(\rho,m)\|_{L^2} \cdot \|\phi\|_{H^1}
\end{equation*}
This norm is finite if we choose $L^2$ norm for $(\rho,m)$ and $H^1$ norm for $\phi$. If we chose $L^2$ norm for $\phi$, the bilinear coupling would have infinite norm. Therefore the corresponding norm of the finite-dimensional coupling on a grid would depend on the grid-size and converge to infinity as grid gets finer. Consequently, the step-sizes that guarantee the convergence of the algorithm would shrink to $0$ and yield an impractically slow algorithm. The $H^1$ norm, on the other hand, yields convergence guarantees and rates that are grid-independent.

For step-sizes $\tau_{\nabla \phi},\tau_{\phi_t},\tau_{\phi(0)},\tau_{\rho},\tau_{m}$, and current iterates $(a^k,\phi^k,\rho^k,m^k,\bar{a}^k,\bar{\phi}^k)$ the update rules for PDHG are
\begin{equation}\label{eq:PDHG_H1}
\begin{cases}
(\rho^{k+1},m^{k+1}) &= \underset{\rho,m}{\mbox{argmax}} ~\mathcal{L} (\bar{\phi}^k,\bar{a}^k,\rho,m)-\frac{1}{2\tau_\rho} \|\rho-\rho^k\|_{L^2_{x,t}}^2-\frac{1}{2\tau_m}\|m-m^k\|_{L^2_{x,t}}^2\\
(a^{k+1},\phi^{k+1})&=\underset{{a,\phi}}{\mbox{argmin}} ~\mathcal{L} (\phi,a,\rho^{k+1},m^{k+1})+\frac{1}{2\tau_{\phi(0)}}\| \phi(\cdot,0)-\phi^k(\cdot,0)\|_{L^2_{x}}^2\\
&+\frac{1}{2\tau_{\nabla \phi}}\|\nabla \phi-\nabla \phi^k\|_{L^2_{x,t}}^2+\frac{1}{2\tau_{\phi_t}}\|\phi_t-\phi^k_t\|_{L^2_{x,t}}^2+\frac{1}{2\tau_a} \|a-a^k\|_{L^2_{t}}^2\\
(\bar{a}^{k+1},\bar{\phi}^{k+1})&=2 (a^{k+1},\phi^{k+1}) - (a^{k},\phi^{k})
\end{cases}
\end{equation}

The critical observation is that the variational problems above are well-posed and easy to solve. In what follows we discuss in details each of the updates.
	
\noindent \textbf{The updates for $(\rho,m)$.} We have that
\begin{equation*}
	\begin{split}
	\frac{\delta \Ll}{\delta \rho}=&-\phi_t-L\left(x,\frac{m}{\rho}\right)+\nabla_v L\left(x,\frac{m}{\rho}\right) \cdot \frac{m}{\rho}-\sum_{i=1}^r a_i(t) f_i(x)\\
	\frac{\delta \Ll}{\delta m}=&-\nabla \phi-\nabla_v L\left(x,\frac{m}{\rho}\right)\\
	\end{split}
\end{equation*}
Therefore, for updating $(\rho,m)$ we must solve the following system
\begin{equation}\label{eq:prox_rho,m}
	\begin{cases}
	\nabla_v L\left(x,\frac{m}{\rho }\right) \cdot \frac{m}{\rho}-L\left(x,\frac{m}{\rho}\right)-\frac{\rho-\rho^k}{\tau_\rho}=\bar{\phi}^k_t+\sum_{i=1}^r \bar{a}^k_i(t) f_i(x)\\
	\nabla_v L\left(x,\frac{m}{\rho}\right)+\frac{m-m^k}{\tau_m}=\nabla \bar{\phi}^k
	\end{cases}
\end{equation}
\begin{remark}\label{rem:prox_parallel}
System \eqref{eq:prox_rho,m} yields decoupled one-dimensional convex optimization problems at the grid-points. Therefore the proximal update for $(\rho,m)$ can de performed efficiently in parallel. This feature is one of the most appealing properties of ADMM \cite{BenamouBrenier2000,bencar'15,bencarsan'17} and PDHG \cite{silva18,silva19,matt19,li18} algorithms for MFG systems and related problems for local couplings. However, direct extensions of aforementioned methods to nonlocal MFG systems do not preserve this property. One of the critical features of our approach is that we preserve this property. We refer to Section \ref{subsec:dim_red} for details.
\end{remark}
For some Lagrangians, \eqref{eq:prox_rho,m} simplifies greatly. For instance, for $L(x,v)=\frac{|v|^2}{2}+Q(x,t)$ we have that $\nabla_v L(x,v)=v$, and \eqref{eq:prox_rho,m} becomes
\begin{equation*}
\begin{cases}
\frac{|m|^2}{2\rho^2}-\frac{\rho-\rho^k}{\tau_\rho}=Q(x,t)+\bar{\phi}^k_t+\sum_{i=1}^r \bar{a}^k_i(t) f_i(x)\\
\frac{m}{\rho}+\frac{m-m^k}{\tau_m}=\nabla \bar{\phi}^k
\end{cases}
\end{equation*}
Furthermore, eliminating $m$ from the second equation, we obtain 
\begin{equation}\label{eq:prox_rho,m_quad}
\begin{cases}
\frac{|m^k+\tau_m \nabla \bar{\phi}^k|^2}{2(\tau_m+\rho)^2}-\frac{\rho-\rho^k}{\tau_\rho}=Q(x,t)+\bar{\phi}^k_t+\sum_{i=1}^r \bar{a}^k_i(t) f_i(x)\\
m=\rho \frac{m^k+\tau_m \nabla \bar{\phi}^k}{\tau_m+\rho}
\end{cases}
\end{equation}
Therefore, we just need to solve a cubic equation for $\rho$ and update $m$ by an explicit formula.
\begin{remark}\label{rem:rho>=0}
As mentioned before, we must add a constraint $\rho \geq 0$ in \eqref{eq:infsupL}. Therefore, equations above must be complemented by the condition $\rho\geq 0$. Furthermore, the expression $L\left(x,\frac{m}{\rho}\right)$ must be understood in the following sense
\begin{equation*}
L\left(x,\frac{m}{\rho}\right)=\begin{cases}
L\left(x,\frac{m}{\rho}\right),~\mbox{when}~\rho>0\\
0,~\mbox{when}~m=0,~\rho=0\\
+\infty,~\mbox{when}~m\neq 0,~\rho=0
\end{cases}
\end{equation*}
The function $\gamma(\rho)=\frac{|m^k+\tau_m \nabla \bar{\phi}^k|^2}{2(\tau_m+\rho)^2}-\frac{\rho-\rho^k}{\tau_\rho},~\rho\geq 0$ is strictly decreasing. Therefore, either $\gamma(0)\geq Q(x,t)+\bar{\phi}^k_t+\sum_{i=1}^r \bar{a}^k_i(t) f_i(x)$ or $\gamma(0)<Q(x,t)+\bar{\phi}^k_t+\sum_{i=1}^r \bar{a}^k_i(t) f_i(x)$.  In the former case, there exists a unique $\rho^{k+1}(x,t)\geq 0$ such that $\gamma(\rho^{k+1}(x,t))=Q(x,t)+\bar{\phi}^k_t(x,t)+\sum_{i=1}^r \bar{a}^k_i(t) f_i(x)$. In the latter case, we set $\rho^{k+1}(x,t)=0$. In both cases, we update $m$ accordingly. 
\end{remark}

\noindent\textbf{The updates for $(a,\phi)$.} We have that
\begin{equation*}
\begin{split}
\frac{\partial \Ll}{\partial a}=& \bK^{-1} a-\left(\int_\Omega f_i(x) \rho(x,t) dx \right)_{i=1}^r\\
\frac{\partial \Ll}{\partial \phi}=& \rho_t+\nabla \cdot m+\left(\rho-\rho_0(x)\right) (dx \rest \Omega)\times \delta_{t=0} -m\cdot \nu~d_{\partial \Omega} x \times dt  
\end{split}
\end{equation*}
where $\nu$ is the outward normal of $\Omega$, and $d_{\partial \Omega}$ is the surface measure of $\partial \Omega$. Note that we only consider variations that preserve the boundary condition $\phi(x,1)=g(x)$. Therefore, to update $(a,\phi)$ we must solve the system
\begin{equation}\label{eq:prox_a,rho}
\begin{cases}
\bK^{-1} a-\left(\int_\Omega f_i(x) \rho^{k+1}(x,t) dx \right)_{i=1}^r+\frac{a-a^k}{\tau_a}=0\\
\rho^{k+1}_t+\nabla \cdot m^{k+1}-\frac{\phi_{tt}-\phi^k_{tt}}{\tau_{\phi_t}}-\frac{\Delta \phi - \Delta \phi^k}{\tau_{\nabla \phi}}=0\\
\rho^{k+1}(x,0)-\rho_0(x)-\frac{\phi_{t}(x,0)-\phi^k_{t}(x,0)}{\tau_{\phi_t}}+\frac{\phi(x,0)-\phi^k(x,0)}{\tau_{\phi(0)}}=0\\
\phi(x,1)=g(x)\\
\left(-m^{k+1}+\frac{\nabla \phi -\nabla \phi^k}{\tau_{\nabla \phi}}\right)\cdot \nu  =0
\end{cases}
\end{equation}
Note that the equations for $a$ and $\phi$ are decoupled. Furthermore, we obtain explicit updates for $a$:
\begin{equation}\label{eq:prox_a}
a^{k+1}=(\tau_a \bK^{-1}+\mathbf{I})^{-1} \left(a^k+\tau_a \left(\int_\Omega f_i(x) \rho^{k+1}(x,t) dx \right)_{i=1}^r\right)
\end{equation}
Finally, to update $\phi$ we must solve the space-time elliptic equation
\begin{equation}\label{eq:prox_phi}
\begin{cases}
\frac{\phi_{tt}}{\tau_{\phi_t}}+\frac{\Delta \phi }{\tau_{\nabla \phi}}=\rho^{k+1}_t+\nabla \cdot m^{k+1}+\frac{\phi^k_{tt}}{\tau_{\phi_t}}+\frac{\Delta \phi^k}{\tau_{\nabla \phi}} \quad \mbox{in}\quad \Omega\times(0,1)\\
\frac{\phi_{t}(x,0)}{\tau_{\phi_t}}-\frac{\phi(x,0)}{\tau_{\phi(0)}}=\rho^{k+1}(x,0)-\rho_0(x)+\frac{\phi^k_{t}(x,0)}{\tau_{\phi_t}}-\frac{\phi^k(x,0)}{\tau_{\phi(0)}} \quad \mbox{in}\quad \Omega\\
\phi(x,1)=g(x)\quad \mbox{in} \quad \Omega\\
\frac{\partial \phi(x,t)}{\partial \nu}=\frac{\partial \phi^k(x,t)}{\partial \nu}+\tau_{\nabla \phi} m^{k+1}\cdot \nu\quad \mbox{in}\quad \partial \Omega \times (0,1)
\end{cases}
\end{equation}
This step can be efficiently performed by the Fast Fourier Transform (FFT).

\subsection{Dimension reduction}\label{subsec:dim_red}

Here, we illustrate the dimension reduction and computational efficiency obtained with our method. Assume that $\Omega$ is bounded, and $K_r,S_r$ are approximations of $K,S$ in the $C^2$ norm. Then we have that
\begin{equation*}
\begin{split}
	\left\|\int_{\Omega} K(\cdot,y) \rho(y) dy- \int_{\Omega} K_r(\cdot,y) \rho(y) dy \right\|_{C^2}\leq \|K-K_r\|_{C^2},\\
	\left\|\int_{\Omega} S(\cdot,y) \rho(y) dy- \int_{\Omega} S_r(\cdot,y) \rho(y) dy \right\|_{C^2}\leq \|S-S_r\|_{C^2},
\end{split}
\end{equation*}
for all $\rho \in \mathcal{P}(\Omega)$. Therefore, if we approximate $f,g$ by $f_r,g_r$ in the $C^2$ norm, we obtain $C^2$ approximations of the terms
\begin{equation*}
f\left(x,\int_\Omega K(x,y)\rho(y,t)dy\right),\quad g\left(x,\int_\Omega S(x,y)\rho(y,1)dy\right)
\end{equation*}
that are uniform in $\rho$. Therefore, $f_r,g_r,K_r,S_r$ produce an approximation of \eqref{eq:main_gen} that is independent of the grid-size or the number of agents.

Furthermore, from the stability theory of \eqref{eq:main_gen} \cite{LasryLions2007,CardaNotes}, we have that solutions of \eqref{eq:main_gen} corresponding to $f_r,g_r,K_r,S_r $ are precompact in $C(\Omega\times[0,1])\times C\left([0,1],\mathcal{P}(\Omega)\right)$, and all accumulation points are solutions corresponding to $f,g,K,S$. Additionally, if the operators
\begin{equation*}
\rho \mapsto f\left(x,\int_\Omega K(x,y)\rho(y)dy\right),\quad \rho \mapsto g\left(x,\int_\Omega S(x,y)\rho(y)dy\right)
\end{equation*}
are monotone, \eqref{eq:main_gen} admits a unique solution, $(\phi,\rho)$, and 
\begin{equation*}
	\begin{split}
	\lim\limits_{r\to \infty}\|\phi_r-\phi\|_{L^\infty}=0,\quad \lim\limits_{r\to \infty} \sup_{t\in [0,1]} W_1(\rho_r(\cdot,t),\rho(\cdot,t))=0,
	\end{split}
\end{equation*}
if $\lim\limits_{r\to \infty}\|\xi-\xi_r\|_{C^2}=0,~ \xi \in \left\{f,g,K,S\right\}$, where $W_1$ is the 1-Wasserstein or Monge-Kantorovich distance in $\mathcal{P}(\Omega)$.

Thus, once we produce approximations of $f,g,K,S$, we obtain an approximation of \eqref{eq:main_gen} that works equally fine across all discretizations. Therefore, once we fix $r$, we can solve the $r$-problem as accurately as we wish without extra cost for fine meshes. Additionally, as we show below, for fixed $r$, the computational cost is on par with those of existing algorithms for local couplings. Of course, the smaller $r$ the better, and the size of $r$ depends on how well $\{f_i\}$ approximate $f,g,K,S$.

We now compare the computational complexity of our method versus direct applications of primal-dual optimization algorithms to solve \eqref{eq:main_gen}. The starting point for these methods is to write \eqref{eq:main} as a convex optimization problem introduced in \cite{LasryLions2007}. More precisely, when $K$ is symmetric, one has that \eqref{eq:main} is equivalent to
\begin{equation}\label{eq:infsupL1_direct}
	\begin{split}
	\inf_{\substack{\phi(x,1) \leq g}} \sup_{\rho\geq 0,m}&  \Bigg\{ -\int_{\Omega}\phi(x,0)\rho_0(x) dx - \int_{\Omega}\int_0^1 \left(\rho \phi_t  + m\cdot \nabla \phi \right) dx dt\\
	&-\int_0^1 \left\{\int_\Omega  \rho L\left(x,\frac{m}{\rho}\right) dx+ \mathcal{F}(\rho(\cdot,t))\right\}dt \Bigg\}\\
	=&\inf_{\substack{\phi(x,1)\leq g\\ \alpha}} \sup_{\rho\geq 0,m}~\mathcal{L}_1(\phi,\rho,m)
	\end{split}	
\end{equation}
where
\begin{equation}\label{eq:Fcal}
\mathcal{F}(\rho)=\frac{1}{2}\int_{\Omega^2} K(x,y)\rho(x)\rho(y)dxdy,\quad \rho \in \mathcal{P}(\Omega).
\end{equation}
One can also work with the convex dual $\mathcal{F}^*$ of $\mathcal{F}$ by introducing a dual variable $\alpha$:
\begin{equation}\label{eq:infsupL2_direct}
\begin{split}
& \inf_{\substack{\phi(x,1) \leq g\\ \alpha}} \sup_{\rho\geq 0,m}   \Bigg\{ \int_0^1 \mathcal{F}^*(\alpha(\cdot,t)) dt -\int_{\Omega}\phi(x,0)\rho_0(x) dx - \int_{\Omega}\int_0^1 \left(\rho \phi_t  + m\cdot \nabla \phi \right) dx dt\\
&-\int_{\Omega}\int_0^1 \rho\left( L\left(x,\frac{m}{\rho}\right) + \alpha \right) dx dt \Bigg\}\\
=&\inf_{\substack{\phi(x,1)\leq g\\ \alpha}} \sup_{\rho\geq 0,m}~\mathcal{L}_2(\phi,\alpha,\rho,m),
\end{split}	
\end{equation}
where $\mathcal{F}^*(\alpha)=\sup_\rho \int_{\Omega} \alpha(x)\rho(x)dx-\mathcal{F}(\rho)$. Therefore, there are two options for solving \eqref{eq:main}: (i) work directly with $\mathcal{F}$ and solve \eqref{eq:infsupL1_direct} or its variants \cite{silva18,silva19}, (ii) work with the dual $\mathcal{F}^*$ and solve \eqref{eq:infsupL2_direct} or its variants \cite{bencar'15,bencarsan'17}. We illustrate that direct applications of both approaches to nonlocal problems lead to computationally expensive updates.

For concreteness, we estimate the computational complexity of the PDHG algorithm proposed here, with and without applying the coefficients method. The analysis of other primal-dual algorithms is analogous. First, we discuss the option of working directly with $\mathcal{F}$ and solving \eqref{eq:infsupL1_direct}. In this case, the proximal update for $(\rho,m)$ would be
\begin{equation*}
\begin{split}
\sup_{\rho\geq 0,m}&  \Bigg\{ -\int_{\Omega}\bar{\phi}^k(x,0)\rho_0(x) dx - \int_{\Omega}\int_0^1 \left(\rho \bar{\phi}^k_t  + m\cdot \nabla \bar{\phi}^k \right) dx dt\\
&-\int_0^1 \left\{\int_\Omega  \rho L\left(x,\frac{m}{\rho}\right) dx+ \mathcal{F}(\rho(\cdot,t))\right\}dt \Bigg\}\\
&-\frac{1}{2\tau_\rho}\|\rho-\rho^k\|_{L^2_{x,t}}^2-\frac{1}{2\tau_m}\|m-m^k\|_{L^2_{x,t}}^2
\end{split}
\end{equation*}
Therefore, we must solve the following system
\begin{equation}\label{eq:prox_rho,m_direct}
\begin{cases}
L\left(x,\frac{m}{\rho}\right)-\nabla_v L\left(x,\frac{m}{\rho }\right) \cdot \frac{m}{\rho}+\frac{\rho-\rho^k}{\tau_\rho}+\delta_\rho \mathcal{F}(\rho)=0\\
\nabla_v L\left(x,\frac{m}{\rho}\right)+\frac{m-m^k}{\tau_m}=\nabla \bar{\phi}^k
\end{cases}
\end{equation}
For local interactions, one has that $\mathcal{F}(\rho)=\int_\Omega F(\rho(x)) dx$ for some $F$, and \eqref{eq:prox_rho,m_direct} becomes
\begin{equation*}
\begin{cases}
L\left(x,\frac{m}{\rho}\right)-\nabla_v L\left(x,\frac{m}{\rho }\right) \cdot \frac{m}{\rho}+\frac{\rho-\rho^k}{\tau_\rho}+F'(\rho)=0\\
\nabla_v L\left(x,\frac{m}{\rho}\right)+\frac{m-m^k}{\tau_m}=\nabla \bar{\phi}^k,
\end{cases}
\end{equation*}
which is a decoupled system of one-dimensional equations that can be solved efficiently at each node. Therefore, the computational complexity of solving this system for local problems is linear in the number of grid-points. However, for nonlocal interactions such as in \eqref{eq:Fcal} we have that $\delta_\rho \mathcal{F}=\int_\Omega K(x,y) \rho(y) dy$, and so \eqref{eq:prox_rho,m_direct} is now a fully coupled (in space) system of nonlinear equations. Additionally, the complexity of the systems grows with the mesh-size. One could approximate the term by $\delta_\rho \mathcal{F}(\rho)$ by $\delta_\rho \mathcal{F}(\rho^k)$ and decouple the system. Nevertheless, this would require a matrix multiplication that yields a superlinear computational cost. Moreover, the proximal step could not be parallelized.

With our method, on the other hand, we obtain fully parallel proximal updates for $(\rho,m)$ \eqref{eq:prox_rho,m} at the expense of solving an $r\times r$ system of equations to update the coefficients $(a_i)$. To assemble the system for $(a_i)$, we need to evaluate terms $\int_\Omega f_i(x) \rho^{k+1}(x,t) dx$ that yields a linear cost in the number of grid-points. Therefore, once $r$ is fixed, we obtain an overall linear cost for updating $(\rho,m)$ and $(a_i)$ which is the case for local interactions. 

Next, we discuss the option of working with $\mathcal{F}^*$ and solving \eqref{eq:infsupL2_direct}. In this case, the proximal update for $\alpha$ would be
\begin{equation*}
	\inf_{\alpha} \int_0^1 \mathcal{F}^*(\alpha(\cdot,t)) dt - \int_{\Omega}\int_0^1 \rho^{k+1} \alpha dx dt +\frac{1}{2\tau_\alpha}\|\alpha-\alpha^k\|^2_{L^2_{x,t}}
\end{equation*}
Therefore, we must solve the following system
\begin{equation}\label{eq:prox_alpha}
	\delta_\alpha \mathcal{F}^*(\alpha)(x,t)-\rho^{k+1}(x,t)+\frac{\alpha(x,t)-\alpha^k(x,t)}{\tau_\alpha}=0
\end{equation}
For local interactions, $\mathcal{F}(\rho)=\int_\Omega F(\rho(x)) dx$, we can calculate $\mathcal{F}^*(\alpha)$ on the continuum level by an explicit formula
\begin{equation*}
	\mathcal{F}^*(\alpha)=\int_\Omega F^*(\alpha(x)) dx,
\end{equation*}
where $F^*$ is the convex dual of $F$. Therefore, $\delta_\alpha \mathcal{F}(\alpha)(x)=(F^*)'(\alpha(x))$, and \eqref{eq:prox_alpha} becomes
\begin{equation*}
(F^*)'(\alpha(x,t))-\rho^{k+1}(x,t)+\frac{\alpha(x,t)-\alpha^k(x,t)}{\tau_\alpha}=0
\end{equation*}
As before, we obtain one-dimensional decoupled equations that can be solved in parallel at grid-points.

In the nonlocal case though, the first issue is that we cannot calculate $\mathcal{F}^*$ analytically on the continuum level unless $K$ is special. However, one can calculate $\mathcal{F}^*$ on the discrete level. Assume that $\left\{x_i\right\}_{i=1}^N$ is some space-discretization. Then we have that
\begin{equation*}
	\mathcal{F}(\rho)=\frac{1}{2} \sum_{i,j} K(x_i,x_j) \rho_{i} \rho_{j},\quad \mathcal{F}^*(\alpha)=\frac{1}{2} \sum_{i,j} Q_{i j} \alpha_{i} \alpha_{j},\quad \delta_\alpha \mathcal{F}^*(\alpha)_{i}= \sum_{j} Q_{ij} \alpha_{j}
\end{equation*}
where $\rho_i=\rho(x_i),~\alpha_i=\alpha(x_i)$, and $Q=(Q_{ij})=(K(x_{i},x_{j}))^{-1}$. Therefore, \eqref{eq:prox_alpha} becomes an $N\times N$ system of linear equations
\begin{equation}\label{eq:prox_alpha_discrete}
	\sum_{j} Q_{ij} \alpha_{j}(t)-\rho^{k+1}_{i} (t)+\frac{\alpha_{i}(t)-\alpha_{i}^k(t)}{\tau_\alpha}=0,\quad i \in \{1,2,\cdots,N \}.
\end{equation}
As before, we obtain a coupled system in the nonlocal case. The solution of this system yields a polynomial computational complexity unless $(K(x_i,x_j))$ is special. For instance, if $(K(x_i,x_j))$ is diagonalizable in a Fourier basis the computational cost is of order $N\log N$ via FFT. 

In our method, on the other hand, we replace dual variables $(\alpha_{i}(t))_{i=1}^N$ by coefficients $(a_i(t))_{i=1}^r$, and \eqref{eq:prox_alpha_discrete} is replaced by an $r\times r$ system \eqref{eq:prox_a}. Therefore, we have to store much less variables and, as mentioned before, obtain a linear computational cost. Additionally, we can calculate the $r\times r$ matrix $(\tau_a \bK^{-1}+{\bf I})^{-1}$ prior to optimization and use it afterward. Moreover, since the size of this matrix is independent of the mesh-size, we do not have to deal with conditioning issues for every mesh-size separately: we can do it once and for all before optimization.

Finally, as we will see in Section \ref{sec:kernels}, a smart choice of basis functions $\{f_i\}$ yields $\bK={\bf I}$. Therefore, the updates for $(a_i(t))_{i=1}^r$ are trivial and there is no need to solve linear systems at all!

\section{Modeling interactions with kernels}\label{sec:kernels}
	
	Here, we discuss modeling aspects of nonlocal MFG systems. In particular, we show how to build kernels to enforce suitable behavior of agents. For that, we draw inspiration from kernel methods in machine learning \cite[Chapter 6]{mohri}.
	
	As mentioned before, \eqref{eq:main} is well posed when $\rho \mapsto \int_\Omega K(x,y)\rho(y)dy$ is monotone. This condition means that agents repel one another and try to minimize the cost
	\begin{equation*}
	\begin{split}
	\phi(x,t)=&\inf_u \int_t^1 \bigg\{L(z(s),u(s))+\int_\Omega K(z(s),y)\rho(y,t)dy\bigg\}ds+g(z(T))\\
	\mbox{s.t.}~&\dot{z}(s)=c(z(s),u(s)),\quad z(t)=x
	\end{split}
	\end{equation*}
	Therefore, $K(x,y)$ is a \textit{similarity measure} between positions $x$ and $y$ that agents try to minimize. Kernel methods in ML study exactly this type of $K$ for data separation. Different choices of $K$ lead to different separations.
	
	The simplest example of $K$ is the inner product, $K(x,y)=x\cdot y$, which is amenable to rigorous mathematical analysis. Natural extensions of the inner product are \textit{positive definite symmetric (PDS)} kernels.
	\begin{definition}
		$K:(x,y)\mapsto K(x,y)$  is a PDS kernel if $(K(x^i,x^j))_{i,j=1}^m$ is symmetric positive semidefinite matrix for all $\{x^i\}_{i=1}^m \subset \mathbb{R}^d$.
	\end{definition}
	Assume $K$ is a continuous PDS. Thus, its symmetric, and for arbitrary $\rho_k=\frac{1}{N}\sum_i w_k^i \delta_{x^i},~k=1,2$ we have that
	\begin{equation*}
	\begin{split}
	&\int_{\Om^2}K(x,y)(\rho_2(x)-\rho_1(x))(\rho_2(y)-\rho_1(y))dxdy\\
	=&\sum_{i,j} K(x^i,x^j)(w_2^i-w_1^i)(w_2^j-w_1^j)\geq 0,
	\end{split}
	\end{equation*}
	and hence $\rho \mapsto \int_\Omega K(x,y)\rho(y)dy$ is monotone.
	
	The discussion above shows that \textit{PDS kernels suit MFG models extremely well}. Thus, we will build various MFG models by choosing suitable PDS kernels. In this context, as we shall see below, the basis $\{f_i\}$ corresponds to feature vectors.
	
	The remarkable fact about PDS kernels is that all of them are inner products. More precisely, $K$ is PDS iff there exists a Hilbert space $\mathcal{H}$ and a mapping $x\mapsto f(x) \in \mathcal{H}$ such that
	\begin{equation}\label{eq:kernel_inner_product}
	K(x,y)=\langle f(x), f(y) \rangle_{\mathcal{H}},\quad \forall x,y
	\end{equation}
	In other words, one can associate points $\{x\}$ in the input space with vectors $\{f(x)\}$ in a Hilbert space so that $K(x,y)$ is precisely the inner product of $f(x)$ and $f(y)$ in $\mathcal{H}$ \cite[Theorem 6.8]{mohri}. The vector $f(x)$ is called the \textit{feature vector} of $x$. If $\mathcal{H}$ is separable we can write $f(x)=(f_1(x),f_2(x),\cdots,f_n(x),\cdots)$ in some basis of $\mathcal{H}$ and $f_1(x),f_2(x),\cdots,f_n(x),\cdots$ will be the \textit{features} of $x$. $\mathcal{H}$ is called a \textit{reproducing kernel Hilbert space (RKHS)}. For $K(x,y)=x\cdot y=\sum x_i y_i $ the features are simply the coordinates.
	
	The RKHS theory blends very well with our coefficients method by providing the basis we need in the form of feature vectors. Indeed, assume that $K$ is a PDS kernel and $\mathcal{H}$ is its RKHS with a basis $\{e_i\}$. Then, we obtain that
	\begin{equation*}
	\begin{split}
	K(x,y)=&\langle f(x), f(y) \rangle_{\mathcal{H}}= \langle \sum_i f_i(x) e_i, \sum_i f_i(y) e_i  \rangle_{\mathcal{H}}\\
	=&\sum_{i,j} \langle e_i,e_j\rangle_{\mathcal{H}} f_i(x) f_j(y)= \sum_{i,j} k_{ij} f_i(x) f_j(y),
	\end{split}
	\end{equation*}
	which is the representation we need. In this case, $\mathbf{K}=(\langle e_i, e_j \rangle_{\mathcal{H}}  )$ is the Gram matrix associated to the basis $\{e_i\}$ in $\mathcal{H}$. Below, we present several common choices for $K$ and provide the basis $\{f_i\}$ and the matrix $\bK$.
	
	\subsection{Maximal spread}\label{subsec:maxspread}
	
	Assume that we want to enforce a maximal spread of the population by penalizing individual agents for being close to the average position of the population. This means that an individual agent faces an optimal control problem 
	\begin{equation*}
	\begin{split}
	&\inf_{\dot{z}=c(z,u)} \int_0^1 \ell(z(s),u(s))-\sum_{i=1}^d \lambda_i \left|z_i(s)-\int_{\Omega} y_i \rho(y,s) dy \right|^2ds+g(z(1))\\
	=&\inf_{\dot{z}=c(z,u)} \int_0^1 \Bigg[\ell (z(s),u(s))-\sum_{i=1}^d \lambda_i |z_i(s)|^2+2\sum_{i=1}^d \lambda_i \left\{\int_{\Omega} z_i(s)\cdot y_i \rho(y,s) dy \right\}\\
	&-\sum_{i=1}^d \lambda_i \left|\int_{\Omega} y_i \rho(y,s) dy \right|^2 \Bigg]ds+g(z(1))
	\end{split}
	\end{equation*}
	Above, $\lambda_1,\lambda_2,\cdots,\lambda_d \geq 0$ signify how much we enforce spreading in each coordinate direction. Additionally, $\ell(x,u)$ is some intrinsic running cost.
	
	Since the term $-\sum_{i=1}^d \lambda_i \left|\int_{\mathbb{R}^d} y_i \rho(y,s) dy \right|^2$ does not depend on the trajectory $z$, the problem above is equivalent to
	\begin{equation*}
	\begin{split}
	\inf_{\dot{z}=c(z,u)}& \int_0^1 \ell (z(s),u(s))-\sum_{i=1}^d \lambda_i |z_i(s)|^2+ \int_{\Omega} \left\{ 2\sum_{i=1}^d \lambda_i z_i(s)\cdot y_i \right\} \rho(y,s) dy dt\\
	&+g(z(1))
	\end{split}
	\end{equation*}
	Therefore, we obtain an MFG system \eqref{eq:main} where
	\begin{equation*}
	\begin{split}
	L(x,u)=&\ell(x,u)-\sum_{i=1}^d \lambda_i |x_i|^2\\
	K(x,y)=&2 \sum_{i=1}^d\lambda_i x_i y_i\\
	H(x,p)=&\sup_u \left\{-p\cdot c(x,u)-L(x,u)\right\}
	\end{split}
	\end{equation*}
	The key point is that $\R^d$ is an RHKS for $K$, and
	\begin{equation*}
	K(x,y)=\sum_{i=1}^d f_i(x) f_i(y)
	\end{equation*}
	where $f_i(x)=\sqrt{2\lambda_i} x_i,~1\leq i \leq d$. Thus, we use these $\{f_i\}$ as the basis in our method. An excellent feature of this choice is that we obtain $\mathbf{K}=\mathbf{K}^{-1}=\mathbf{I}$ which yields a trivial update rule \eqref{eq:prox_a} for $a$ that reads
\begin{equation*}
a^{k+1}_i(t) = \frac{\tau_{a} \int_{\Omega} f_i(x) \rho^{k+1}(x,t) dx + a_i^{k}(t)}{\tau_a +1}
\end{equation*}	
	
	\subsection{Gaussian repulsion}\label{subsec:gauss}
	
	Another common choice for PDS kernels are Gaussians
	\begin{equation*}
	K(x,y)=\mu \prod_{i=1}^d\exp \left(-\frac{|x_i-y_i|^2}{2\sigma_i^2}\right),
	\end{equation*}
	for some $\mu,\sigma_1,\sigma_2,\cdots,\sigma_d >0$. The parameter $\sigma_i$ signifies how repulsive are the agents in $i$-th coordinate direction.
	
	As before, we will try to find a suitable expansion of $K$. Using the power series expansion of $e^x$ one can show that 
	\begin{equation*}
	\begin{split}
	K(x,y)=&\sum_{\alpha_1,\alpha_2,\cdots,\alpha_d \geq 0}  \left\{\sqrt{\mu} e^{-\sum_{i=1}^d\frac{|x_i|^2}{2\sigma_i^2}} \prod_{i=1}^{d}\frac{x_i^{\alpha_i}}{\sigma_i^{\alpha_i}\alpha_i!}\right\} \cdot \left\{\sqrt{\mu}  e^{-\sum_{i=1}^d\frac{|y_i|^2}{2\sigma_i^2}}\prod_{i=1}^{d}\frac{y_i^{\alpha_i}}{\sigma_i^{\alpha_i}\alpha_i!}\right\}\\
	=&\sum_{\alpha_1,\alpha_2,\cdots,\alpha_d \geq 0} f_{\alpha_1,\alpha_2,\cdots,\alpha_d}(x) f_{\alpha_1,\alpha_2,\cdots,\alpha_d}(y) ,
	\end{split}
	\end{equation*}
	where
	\begin{equation*}
	f_{\alpha_1,\alpha_2,\cdots,\alpha_d}(x)= \sqrt{\mu} e^{-\sum_{i=1}^d\frac{|x_i|^2}{2\sigma_i^2}} \prod_{i=1}^{d}\frac{x_i^{\alpha_i}}{\sigma_i^{\alpha_i}\alpha_i!},\quad \alpha_1,\alpha_2,\cdots,\alpha_d \geq 0.
	\end{equation*}
	Hence, we choose $\{f_{\alpha_1,\alpha_2,\cdots,\alpha_d}\}$ as the basis for our coefficients method and choose $n$ to approximate $K$ with functions of order $\sum_{i=1}^d\alpha_i\leq n$. As before, an excellent feature of this basis is that $\mathbf{K}=\mathbf{K}^{-1}=\mathbf{I}$.
	
	\subsection*{Interactions in sub-regions}\label{subsec:subregions}
	
	Methods described above also provide flexible framework to model interactions within sub-regions. Assume that $\Omega_1,\Omega_2$ are complementary in $\Omega$. Furthermore, assume that agents in $\Omega_i$ interact only with those in $\Omega_i$ for $i=1,2$. There is a straightforward way of extending the framework above to this setting.
	
	Suppose that kernels modeling the interaction in $\Omega_1,\Omega_2$ are $K_1,K_2$, respectively. Additionally, assume that the basis for $K_1$ is $\{f^1_i\}$, and the one for $K_2$ is $\{f^2_i\}$ that can be of a different size. Thus, we want to construct $K$ such that
	\begin{equation*}
	K(x,y)=\begin{cases}
	K_i(x,y),\quad (x,y)\in \Omega_i\times \Omega_i\\
	0,\quad \mbox{otherwise} 
	\end{cases}
	\end{equation*}
	Furthermore, we want to construct a basis for $K$ out of $\{f^1_i\}$ and $\{f^2_i\}$. These can be done as follows:
	\begin{equation*}
	\begin{split}
	K(x,y)=&K_1(x,y)\cdot \chi_{\Omega_1}(x)\chi_{\Omega_1}(y)+K_2(x,y)\cdot \chi_{\Omega_2}(x)\chi_{\Omega_2}(y)=\\
	=&\sum_{ij} k^1_{ij} f^1_i(x) f^1_j(y) \chi_{\Omega_1}(x)\chi_{\Omega_1}(y)+\sum_{ij} k^2_{ij} f^2_i(x) f^1_j(y) \chi_{\Omega_2}(x)\chi_{\Omega_2}(y)\\
	=&\sum_{ij} k^1_{ij} f^1_i(x) \chi_{\Omega_1}(x) f^1_j(y) \chi_{\Omega_1}(y)+\sum_{ij} k^2_{ij} f^2_i(x) \chi_{\Omega_2}(x) f^2_j(y) \chi_{\Omega_2}(y),
	\end{split}
	\end{equation*}
	where $\chi_A$ is the characteristic function of $A$. Therefore, the basis for $K$ can be obtained from the ones of $K_1,K_2$ by simply restricting them into subdomains and combining:
	\begin{equation*}
	\big\{  f^1_i(x) \chi_{\Omega_1}(x),~ f^2_i(x) \chi_{\Omega_2}(x) \big\}
	\end{equation*}
	Furthermore, we have that
	\begin{equation*}
	\mathbf{K}=\begin{pmatrix}
	\mathbf{K}_1& \mathbf{0}\\
	\mathbf{0} & \mathbf{K}_2
	\end{pmatrix},\quad \mathbf{K}^{-1}=\begin{pmatrix}
	\mathbf{K}_1^{-1}& \mathbf{0}\\
	\mathbf{0} & \mathbf{K}_2^{-1}
	\end{pmatrix},
	\end{equation*}
	where $\bK_i=(k^i_{ij}),~i=1,2$.  Therefore, low complexity matrices for $K_1,K_2$ yield a low complexity matrix for $K$.
	
	In case we have multiple regions $\Omega_1,\Omega_2,\cdots, \Omega_N$ with kernels $K_1,K_2,\cdots,K_N$ and bases $\{f^1_i\},\{f^2_i\},\cdots \{f^N_i\}$  we obtain a basis
	\begin{equation*}
	\big\{  f^1_i(x) \chi_{\Omega_1}(x),~ f^2_i(x) \chi_{\Omega_2}(x),~\cdots, f^N_i(x) \chi_{\Omega_N}(x) \big\},
	\end{equation*}
	and
	\begin{equation*}
	\bK=\begin{pmatrix}
	\bK_1& \mathbf{0} & \cdots & \mathbf{0}\\
	\mathbf{0} & \bK_2 & \cdots & \mathbf{0}\\
	\vdots & \vdots & \vdots & \vdots\\
	\mathbf{0} & \cdots & \bK_{N-1} & \mathbf{0} \\
	\mathbf{0} & \cdots & \mathbf{0} & \bK_N \\
	\end{pmatrix},\quad \bK^{-1}=\begin{pmatrix}
	\bK_1^{-1}& \mathbf{0} & \cdots & \mathbf{0}\\
	\mathbf{0} & \bK_2^{-1} & \cdots & \mathbf{0}\\
	\vdots & \vdots & \vdots & \vdots\\
	\mathbf{0} & \cdots & \bK_{N-1}^{-1} & \mathbf{0} \\
	\mathbf{0} & \cdots & \mathbf{0} & \bK_N^{-1}\\
	\end{pmatrix},
	\end{equation*}
	where $\bK_i$ is the matrix corresponding to $K_i$.

	\subsection{Interactions given by differential operators}\label{subsec:diffop}
	
	Finally, we demonstrate how our methods work for interactions given by differential operators. In a seminal paper on numerical methods for MFG \cite{achdou10}, the authors consider an interaction term
	\begin{equation}\label{eq:V}
	V[\rho]=\mu (I-\Delta)^{-2} \rho,
	\end{equation}
	where $\mu>0$,  $\Delta$ is the Laplacian operator, and the problem is set on a flat torus $\mathbb{T}^d$. We have that
	\begin{equation*}
	V[\rho]=\int_{\mathbb{T}^d} \Gamma(x-y)\rho(y)dy,
	\end{equation*}
	where $\Gamma$ is the fundamental solution; that is,
	\begin{equation}\label{eq:fundamental}
	(I-\Delta)^2 \Gamma=\mu \delta_0.
	\end{equation}
	Thus, we have that $K(x,y)=\Gamma(x-y)$. As pointed out in \cite{nursaude18}, for convolutions on a torus, the appropriate basis is the trigonometric one. Thus, we need to expand $\Gamma$ into Fourier series with respect to functions $\{\cos(2\pi \alpha\cdot x),\sin(2\pi \alpha \cdot x) \}$ where $\alpha=(\alpha_1,\alpha_2,\cdots,\alpha_d)$. Furthermore, by the symmetry we have that $\Gamma(x)=\Gamma(-x)$. Therefore, the expansion of $\Gamma$ contains only even functions; that is,
	\begin{equation*}
	\Gamma(x)=\sum_{\alpha\geq 0} \gamma_\alpha \cos(2\pi \alpha\cdot x)
	\end{equation*}
	Next, solving \eqref{eq:fundamental} in a Fourier space yields
	\begin{equation*}
	\gamma_0=\mu,\quad \gamma_\alpha=\frac{2\mu}{1+8\pi^2 |\alpha|^2+16 \pi^4 |\alpha|^4},\quad \alpha>0,
	\end{equation*}
	where $|\alpha|^2=\sum_{i=1}^d \alpha_i^2$. Furthermore, we have that
	\begin{equation*}
	\begin{split}
	K(x,y)=&\Gamma(x-y)=\sum_{\alpha\geq 0} \gamma_\alpha \cos(2\pi \alpha \cdot (x-y) )\\
	=& \sum_{\alpha\geq 0} \bigg(\gamma_\alpha \cos(2\pi \alpha \cdot x ) \cos(2\pi \alpha \cdot y )+\gamma_\alpha \sin(2\pi \alpha \cdot x ) \sin(2\pi \alpha \cdot y ) \bigg)\\
	=&\sum_{\alpha \geq 0} f^{\mathrm{cos}}_\alpha(x) f^{\mathrm{cos}}_\alpha(y) +\sum_{\alpha > 0} f^{\mathrm{sin}}_\alpha(x) f^{\mathrm{sin}}_\alpha(y),
	\end{split}
	\end{equation*}
	where
	\begin{equation*}
	f^{\mathrm{cos}}_\alpha(x)=\sqrt{\gamma_\alpha} \cos(2\pi \alpha \cdot x ),\quad f^{\mathrm{sin}}_\alpha(x)=\sqrt{\gamma_\alpha} \sin(2\pi \alpha \cdot x ),\quad \alpha\geq 0.
	\end{equation*}
	Therefore, we choose $\{f^{\mathrm{cos}}_{\alpha},~f^{\mathrm{sin}}_{\alpha}\}$ as the basis for our coefficients method and choose $n$ to approximate $K$ with functions of order $\sum_{i=1}^d \alpha_i \leq n$. Again, this choice renders $\mathbf{K}=\mathbf{K}^{-1}=\mathbf{I}$.

\section{Potential applications to multi-agent trajectory planning problems}\label{sec:multi_agent_control}

Here, we discuss potential applications of our methods to multi-agent trajectory planning problems. We start by a brief derivation of \eqref{eq:main_gen} and \eqref{eq:main}. Assume that we have a swarm of agents where agent $i\in \{1,2,\cdots,N\}$ aims at minimizing a cost
\begin{equation}\label{eq:N_oc}
\begin{split}
&\inf_{u_i} \int_t^1 L(z_i(s),u_i(s),s)+f_i(z_i(s),z_{-i}(s))ds+g_i(z_i(1),z_{-i}(1))\\
\mbox{s.t.}~& \dot{z}_i(s)=c(z_i(s),u_i(s)),\quad z_i(t)=x_i
\end{split}
\end{equation}
Above, $z_{-i}=(z_j)_{j\neq i}$, and $f,g$ model interactions between the agents. This problem leads to a system of $N$ coupled HJBs that is extremely challenging to solve, especially in high-dimensions and for many agents.

The MFG framework provides a solution to this problem by assuming symmetric interactions and considering the continuum limit $N=\infty$. More precisely, if we suppose that
\begin{equation*}
f_i(z_i,z_{-i})=f\left(z_i,\frac{1}{N-1}\sum_{j\neq i}\delta_{z_j}\right),\quad 	g_i(z_i,z_{-i})=g\left(z_i,\frac{1}{N-1}\sum_{j\neq i}\delta_{z_j}\right),
\end{equation*}
and formally pass to the limit when $N\to \infty$ we obtain a system where a generic agent solves an optimal control problem
\begin{equation}\label{eq:indiv_oc}
\begin{split}
\phi(x,t)=&\inf_{u} \int_t^1 L(z(s),u(s),s)+f(z(s),\rho(\cdot,s))ds+g(z(1),\rho(\cdot,1))\\
\mbox{s.t.}~& \dot{z}(s)=c(z(s),u(s)),\quad z(t)=x,
\end{split}
\end{equation}
where $\rho(\cdot,s)$ is the distribution of population at time $s$. We have that $\phi$ solves the HJB equation
\begin{equation*}
\begin{cases}
-\phi_t+\sup_u \left\{-\nabla \phi(x,t) \cdot c(x,u)-L(x,u,t)\right\}=f(x,\rho(x,t))\\
\phi(x,1)=g(x,\rho(x,1)),
\end{cases}
\end{equation*}  
and the optimal control, $u^*$, is given by the Pontryagin Maximum Principle:
\begin{equation*}
u^*(x,t)\in \mbox{argmax}_u \left\{-\nabla\phi(x,t)\cdot c(x,u)-L(x,u,t) \right\}
\end{equation*}
Furthermore, $\rho$ satisfies the continuity equation
\begin{equation*}
\begin{cases}
\rho_t(x,t)+\nabla \cdot \left(\rho(x,t) c(x,u^*(x,t))\right)=0\\
\rho(x,0)=\rho_0(x),
\end{cases}
\end{equation*}
where $\rho_0$ is the initial distribution of the agents. Collecting all equations together we obtain the MFG system
\begin{equation}\label{eq:main_oc}
\begin{cases}
-\phi_t+\sup_u \left\{-\nabla \phi \cdot c(x,u)-L(x,u,t)\right\}=f(x,\rho(x,t))\\
\rho_t(x,t)+\nabla \cdot \left(\rho(x,t) c(x,u^*(x,t))\right)=0\\
u^*(x,t)\in \mbox{argmax}_u \left\{-\nabla\phi(x,t)\cdot c(x,u)-L(x,u) \right\}\\
\rho(x,0)=\rho_0(x),~\phi(x,1)=g(x,\rho(x,1))
\end{cases}
\end{equation}
\begin{remark}
	Equations \eqref{eq:main_gen}, \eqref{eq:main} correspond to the case $c(x,u)=u$, and $u\mapsto L(x,u)$ convex, for which we have a rigorous mathematical analysis \cite{LasryLions2007}. 
\end{remark}
The appealing feature of MFG systems is that instead of solving a highly coupled system of $N$ HJBs in $d\times N$ dimensions we have to solve a single HJB coupled with a continuity equation in $d$ dimensions. More importantly, the MFG optimal control, $u^*$, yields \textit{near optimal} controls for the $N$ agent problem \eqref{eq:N_oc}.

Of course, the performance of the MFG control in \eqref{eq:N_oc} depends on $N$ and gets better as $N$ grows. Nevertheless, it still makes sense to apply MFG controls because they are faster to generate and provide appealing properties such as no-collision trajectories. For instance, if $c(x,u)=u$, and $L(x,u)=\frac{|u^2|}{2}+Q(x,t)$ for some smooth $Q$, one can show that trajectories corresponding to $u^*(x,t)$ do not intersect \cite[Lemma 4.13]{CardaNotes}.

In the context above, our method may provide a flexible way of augmenting existing solution methods for single-agent trajectory planning problems to generate MFG optimal controls for multi-agent problems.

Indeed, Theorem \ref{thm:optimization} asserts that, under the settings of this paper, \eqref{eq:main_oc} is equivalent to the optimization problem
\begin{equation}\label{eq:optimization_oc}
\begin{split}
&\inf_{a\in C\left([0,1];\mathbb{R}^r\right)}  \frac{\langle \mathbf{K}^{-1} a, a \rangle }{2}- \int_\Omega \phi_a(x,0)\rho_0(x)dx
\end{split}
\end{equation}
where $\phi_a$ solves the HJB
\begin{equation}\label{eq:HJB_oc}
\begin{cases}
-\phi_t+\sup_u \left\{-\nabla \phi(x,t) \cdot c(x,u)-L(x,u,t)\right\}=\sum_{i=1}^r a_i(t)f_i(x)\\
\phi(x,1)=g(x),
\end{cases}
\end{equation}
In Section \ref{sec:pdhg}, we showed how to apply a PDHG algorithm to solve \eqref{eq:optimization_oc}. Here, we argue that virtually any HJB solver (single-agent trajectory-planner) can be augmented to solve \eqref{eq:optimization_oc}.

We propose solving \eqref{eq:optimization_oc} by some type of gradient descent on $a=(a_i)$. For that, we fix an iterate $a^{current}$ and run any single-agent trajectory planning algorithm to solve \eqref{eq:HJB_oc} for $a=a^{current}$ and generate an optimal control $u^{current}$. Then we solve the forward continuity equation
\begin{equation*}
\begin{cases}
\rho_t(x,t)+\nabla \cdot \left(\rho(x,t) c(x,u^{current}(x,t))\right)=0\\
\rho(x,0)=\rho_0(x),
\end{cases}
\end{equation*}
and generate $\rho^{current}$. Finally, we update $a=(a_i)$ by a gradient descent step using \eqref{eq:dGa}:
\begin{equation*}
a^{new}_i(t)=a^{current}_i(t)-h \sum_{j=1}^r k_{ij} a_j^{current}(t)+h \int_{\Omega} f_i(x)\rho^{current}(x,t)dx,
\end{equation*}
where $h>0$ is the descent step-sizes.
\begin{remark}
	Several remarks are in order.
	\begin{enumerate}
		\item Explicit gradient descent steps and exact solutions $u^{current}$ can be replaced by implicit (proximal) steps and approximate solutions as in the PDHG here and \cite{nursaude18}.
		\item The approach above works for both Eulerian and Lagrangian solvers. For latter, the terms $\int_\Omega f_i(x)\rho^{current}(x,t)$ are simply averages of $f_i$-s on the trajectories of particles \cite{nursaude18}.
		\item The number of coefficients, $r$, does not depend on the number of agents, and we always need to (approximately) solve \textit{one decoupled} HJB at each iteration. 
	\end{enumerate}
\end{remark}

Finally, we observe that the methods discussed here also work the other way around; that is, optimal control solvers can be easily adapted to solve MFG problems.

\section{Numerical experiments}\label{sec:numerical_examples}
	
Here, we present several numerical experiments in a two-dimensional case for kernels and bases discussed in Section \ref{sec:kernels}. We take $\Omega\times[0,T] = [0,1]^2 \times [0,1]$ and choose a uniform space-grid with $N_x=64$ points per dimension and a uniform time-grid with $N_t=32$ points in all our examples.

Given $N_{x1},N_{x2},N_t$, we have $\Delta x_1 = \frac{1}{N_{x1}}$, $\Delta x_2 = \frac{1}{N_{x2}}$, $\Delta t = \frac{1}{N_{t}}$.
For $x_1 = i \Delta x_1, x_2 = j \Delta x_2,t_l =l \Delta t $, define  
\begin{equation*}
\begin{split}
\rho_{i,j}^l &= \rho(x_i,x_j,t_l) ~\quad 1\leq i \leq N_{x1}, 1\leq j\leq N_{x2}, 1\leq l \leq N_t\\
m_{1,i+\frac{1}{2},j+\frac{1}{2}}^l &= m_{x1}(x_{i+\frac{1}{2}},x_{j+\frac{1}{2}},t_l)  ~\quad 1\leq i \leq N_{x1}, 1\leq j\leq N_{x2}, 1\leq l \leq N_t\\
m_{2,i+\frac{1}{2},j+\frac{1}{2}}^l &= m_{x2}(x_{i+\frac{1}{2}},x_{j+\frac{1}{2}},t_l)  ~\quad 1\leq i \leq N_{x1}, 1\leq j\leq N_{x2}, 1\leq l \leq N_t\\
\phi_{i,j}^l &= \phi(x_i,x_j, t_l)  ~\quad 1\leq i \leq N_{x1}, 1\leq j\leq N_{x2}, 1\leq l \leq N_t\\
a_{k,l}  &= a_{i}(t_l) ~\quad 1\leq l \leq N_t, 1\leq k \leq r\\
f_{k,i,j} &= f_k(x_i,x_j) ~\quad 1\leq i \leq N_{x1}, 1\leq j\leq N_{x2}, 1\leq k \leq r\\
\rho^0_{i,j} &= \rho_0(x_i,x_j) ~\quad 1\leq i \leq N_{x1}, 1\leq j\leq N_{x2} \\
g_{i,j} &= g(x_i,x_j) ~\quad 1\leq i \leq N_{x1}, 1\leq j\leq N_{x2}
\end{split}
\end{equation*}
To satisfy the Neumann boundary condition, we have
\begin{equation*}
\begin{split}
m_{1,i+\frac{1}{2},j+\frac{1}{2}}^l &= 0  ~\quad i=1,N_{x1}, 1\leq j\leq N_{x2}, 1\leq l \leq N_t\\	  m_{2,i+\frac{1}{2},j+\frac{1}{2}}^l &= 0  ~\quad 1\leq i \leq N_{x1},j=1, N_{x2}, 1\leq l \leq N_t\\
\end{split}
\end{equation*}
The Fokker-Planck equation discretized with forward difference in time as follows:
\begin{equation*}
\frac{1}{\Delta t} \left( \rho_{i,j}^{l+1} - \rho_{i,j}^{l}\right) + \frac{1}{\Delta x_1} \left( m_{i+\frac{1}{2},j}^{l} -m_{i-\frac{1}{2},j}^{l}  \right) + \frac{1}{\Delta x_2} \left( m_{i, j+\frac{1}{2}}^{l} -m_{i,j-\frac{1}{2}}^{l}  \right)= 0,
\end{equation*}
The HJB equation is discretized with backward difference in time as follows:
\begin{equation*}
\begin{split}
-\frac{1}{\Delta t} \left( \phi_j^{l} - \phi_j^{l-1}\right) + H\left(x_j, \nabla_x \phi_{i,j}^{l}  \right) = \sum_{k=1}^r a_{k,l} f_{k,i,j},\\
\text{with}  \quad \nabla_x \phi_{i,j}^{l}  = \left[\frac{\phi_{i+1,j}^{l} -\phi_{i,j}^{l}}{\Delta x_1},\frac{\phi_{i,j+1}^{l} -\phi_{i,j}^{l}}{\Delta {x_2}} \right]^T.
\end{split}
\end{equation*}
With above finite difference scheme, we are ready to define the objective function of our min-max problem $\mathcal{L}(\phi,a,\rho,m)$ in discrete form:
\begin{equation*}
\begin{split}
&{\mathcal{L}}(\phi,a,\rho,m) = \frac{\Delta t}{2}\sum_l \sum_{k_1,k_2=1}^r a_{k_1,l}a_{k_2,l}r_{k_1 k_2} -\Delta x_1 \Delta x_2  \sum_{i,j}\phi_{i,j,1}\rho^0_{i,j} \\&- \Delta x_1 \Delta x_2 {\Delta t} \sum_{i,j} \sum_l \Bigg(\left(\rho_{i,j}^l \frac{\left( \phi_j^{l} - \phi_j^{l-1}\right)}{\Delta t} +  m_{i+\frac{1}{2},j+\frac{1}{2}}^l \cdot \nabla_x \phi_{i,j}^{l}  \right)\\
& + \rho_{i,j}^l \left(L(x,\frac{ m_{i+\frac{1}{2},j+\frac{1}{2}}^l}{\rho_{i,j}^l})+\sum_{k=1}^r a_{k,l} f_{k,i,j}\right)\Bigg),\\
&\quad \text{with }\quad m_{i,j}^l = \left[ m_{1,i+\frac{1}{2},j+\frac{1}{2}}^l, m_{2,i+\frac{1}{2},j+\frac{1}{2}}^l \right]^T,\\
\end{split}
\end{equation*}
where $(r_{k_1 k_2})=\bK^{-1}$.

\subsection{Maximal spread}
We consider a maximal spread model on from Section \ref{subsec:maxspread} in the domain $\Omega \times [0,T] = [0,1]^2 \times [0,1]$. We denote by $\rho_G(c_1,c_2,\sigma_G)$ the density of a homogeneous normal distribution centered at $(c_1,c_2)$ with variance $\sigma_G^2$. We set the initial-terminal conditions for our MFG system to be
	\begin{equation*}
	\begin{split}
	\rho_0(x_1, x_2) & = \rho_G(0.5,0.9,0.2)\\
	%\frac{1}{2 \pi \sigma^2 } \exp \left(-\dfrac{(x_1 -0.5)^2 + (x_2-0.9)^2 }{2\sigma^2}\right), \quad \sigma = 0.2\\
	g(x_1, x_2) & = 2\exp\left(-10 \left(x_1-0.5 \right)^2- 
		\left(x_2-0.1\right)^2\right) \left(\left(x_2-0.1 \right)^2 -1\right).\\
	\end{split}
	\end{equation*}	
Furthermore, we set 
	\begin{equation*}
		\begin{split}
		L(x,v) &= \frac{1}{2}\|v\|^2 + 10^3 \left(\max \left(|x_1-0.5|,|x_2-0.5| \right)\right)^8\\
		%Q(x) &= 10^3 \left(\max \left(|x_1-0.5|,|x_2-0.5| \right)\right)^8.
		\end{split}
	\end{equation*} 
%& Q(x) =c \left(\max(|x_1-0.5|,|x_2-0.5|)\right)^8,\quad c=5e3 \\
\begin{figure} [!htbp]
	\centering
	\begin{subfigure}{\linewidth}
		\includegraphics[width=1.0\textwidth,trim=0 0 0 0, clip=false]{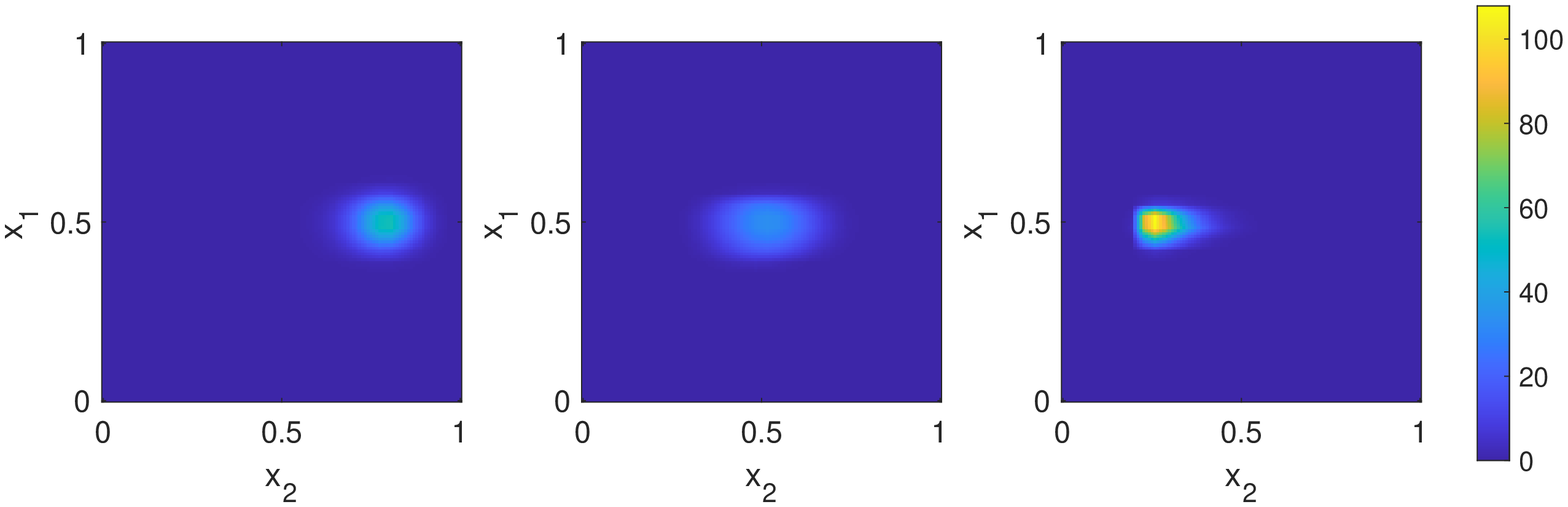}
	\vspace{-0.7cm}
	\caption{$\lambda_1 = 0.1, \lambda_2= 0.1 $}
	\end{subfigure}
	\begin{subfigure}{\linewidth}
		\includegraphics[width=1.0\textwidth,trim=0 0 0 0, clip=false]{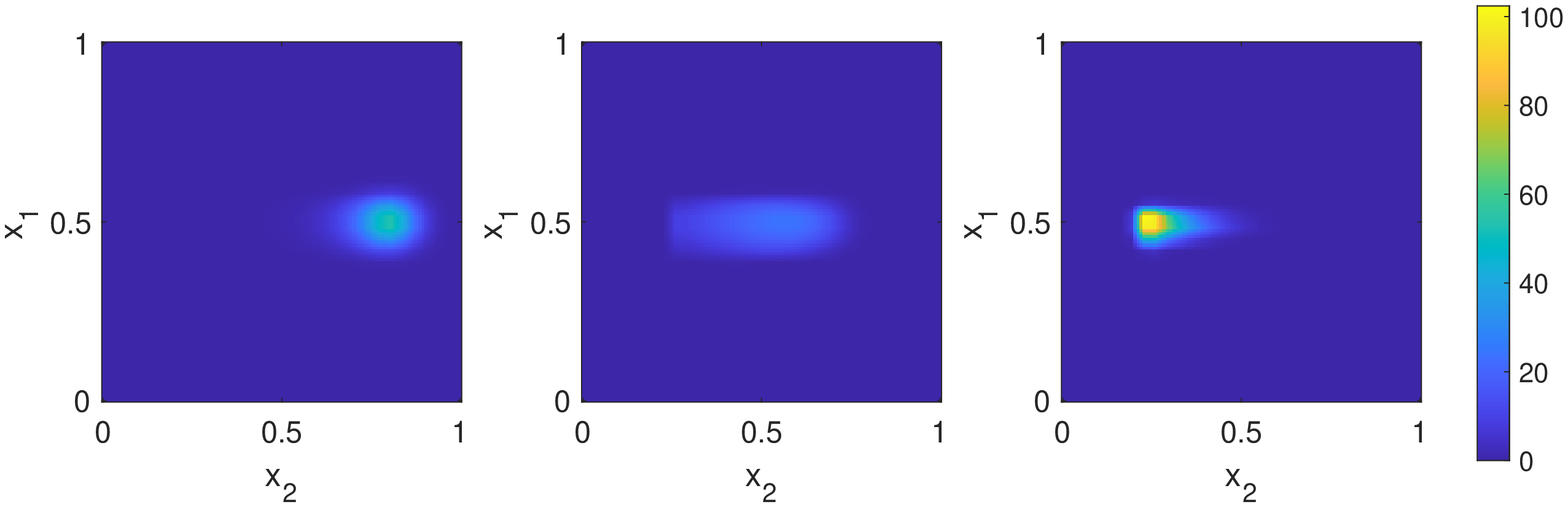}
	\vspace{-0.7cm}
	\caption{$\lambda_1 = 0.1, \lambda_2= 4 $}
\end{subfigure}
		\begin{subfigure}{\linewidth}
		\includegraphics[width=1.0\textwidth,trim=0 0 0 0, clip=false]{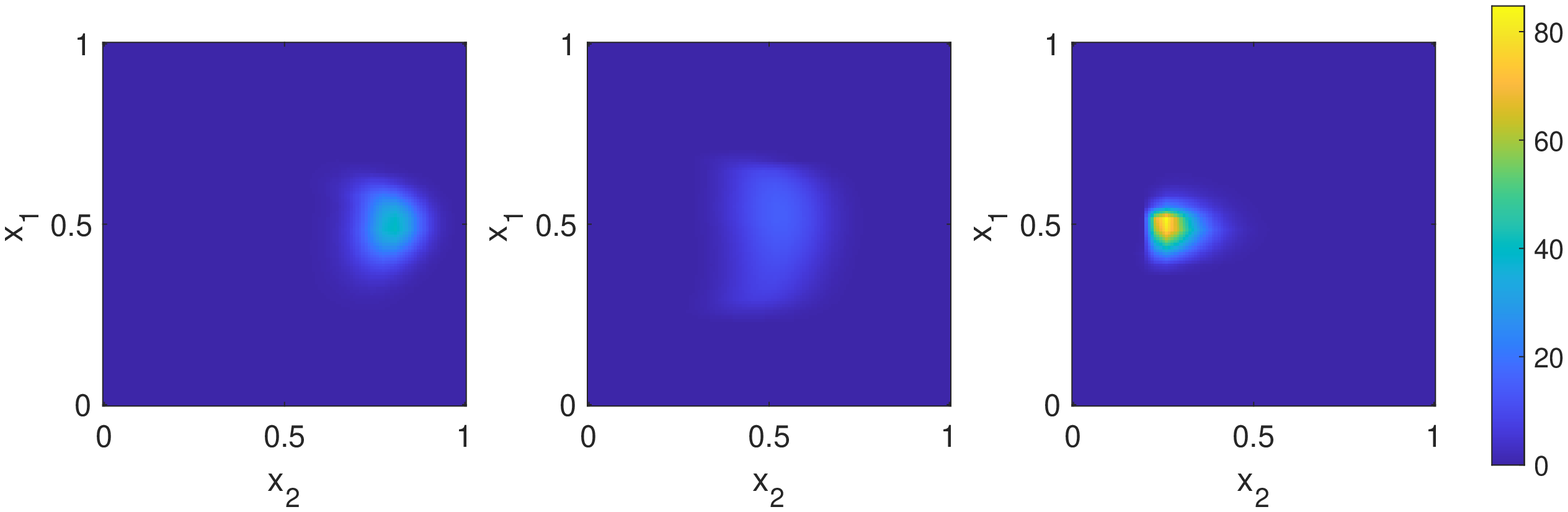}
	\vspace{-0.7cm}
	\caption{$\lambda_1 = 4, \lambda_2= 0.1 $}
\end{subfigure}
		\begin{subfigure}{\linewidth}
		\includegraphics[width=1.0\textwidth,trim=0 0 0 0, clip=false]{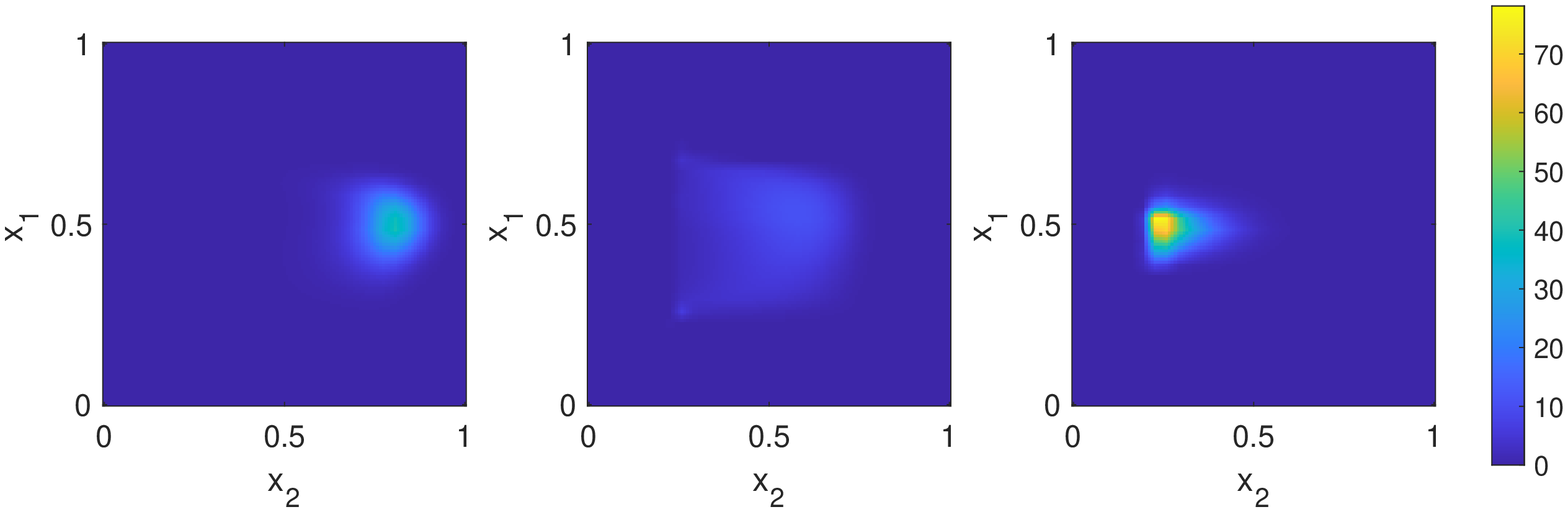}
	\vspace{-0.7cm}
	\caption{$\lambda_1 = 4, \lambda_2= 4 $}
\end{subfigure}
\caption{\textbf{Maximal spread.} MFG solution $\rho(x,t)$ at $t=0.1, 0.5, 0.9$ with different choice of $\lambda_i$}
	\label{fig:maximal_spread}
\end{figure}
%    \begin{figure} [!htbp]
%	\centering
%	\includegraphics[width=1.0\textwidth,trim=0 0 0 0, clip=false]{maximal_spread_eg_4_lambdas.eps}
%	\vspace{-0.7cm}
%	\caption{MFG solution $\rho(x,0.5)$ for different choice of $\lambda_i$}
%	\label{fig:maximal_spread}
%\end{figure}
We have computed the MFG solutions for four choices of parameters
\begin{equation*}
(\lambda_{1}, \lambda_2)\in \left\{ (0.1,0.1),(0.1,4),(4,0.1) ,(4,4) \right\}
\end{equation*}
The results are shown in Figure \ref{fig:maximal_spread}. We can see that, in accordance to theory, larger $\lambda_i$ prompt larger spread in $x_i$ directions. Additionally, we see the flexibility of our method for modeling interactions that are heterogeneous in different directions.

 %the MFG system has $64^2 \times 32$ in space-time. 
	
	\subsection{Gaussian repulsion with static obstacles}\label{subsec:gaussian_static}
\begin{figure} [!htbp]
		\centering
		\begin{subfigure}{\linewidth}
			\includegraphics[width=1.0\textwidth,trim=0 0 0 0, clip=false]{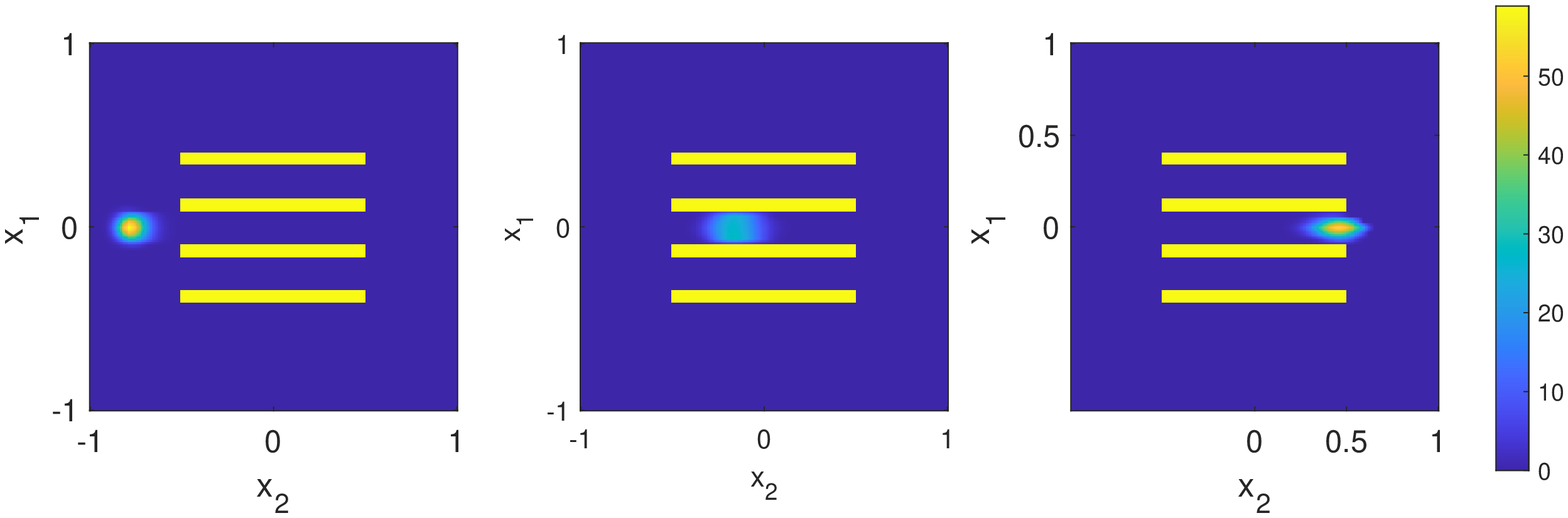}
			\vspace{-0.7cm}
			\caption{$\sigma_1 = 0.8, \sigma_2 = 0.8, \mu = 0.1 $}
		\end{subfigure}
		\begin{subfigure}{\linewidth}
			\includegraphics[width=1.0\textwidth,trim=0 0 0 0, clip=false]{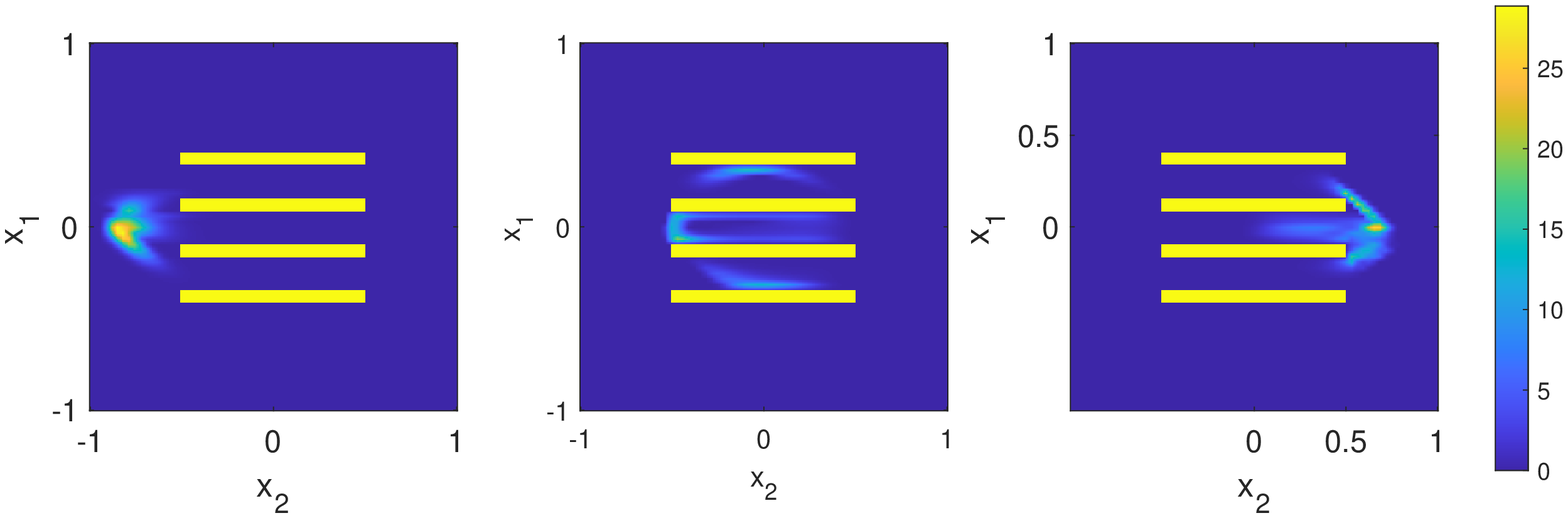}
			\vspace{-0.7cm}
			\caption{$\sigma_1 = 0.2, \sigma_2 = 0.2, \mu = 5$}
		\end{subfigure}
		\begin{subfigure}{\linewidth}
			\includegraphics[width=1.0\textwidth,trim=0 0 0 0, clip=false]{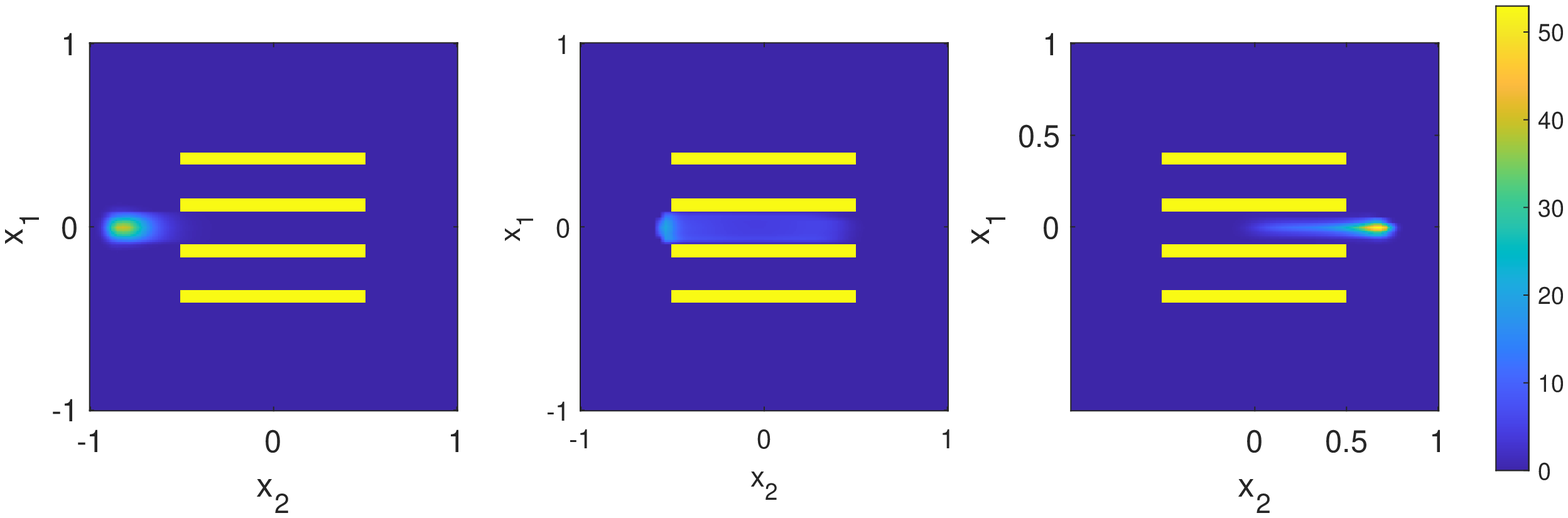}
			\vspace{-0.7cm}
			\caption{$\sigma_1 =0.5, \sigma_2 = 0.2, \mu = 5 $}
		\end{subfigure}
		\begin{subfigure}{\linewidth}
			\includegraphics[width=1.0\textwidth,trim=0 0 0 0, clip=false]{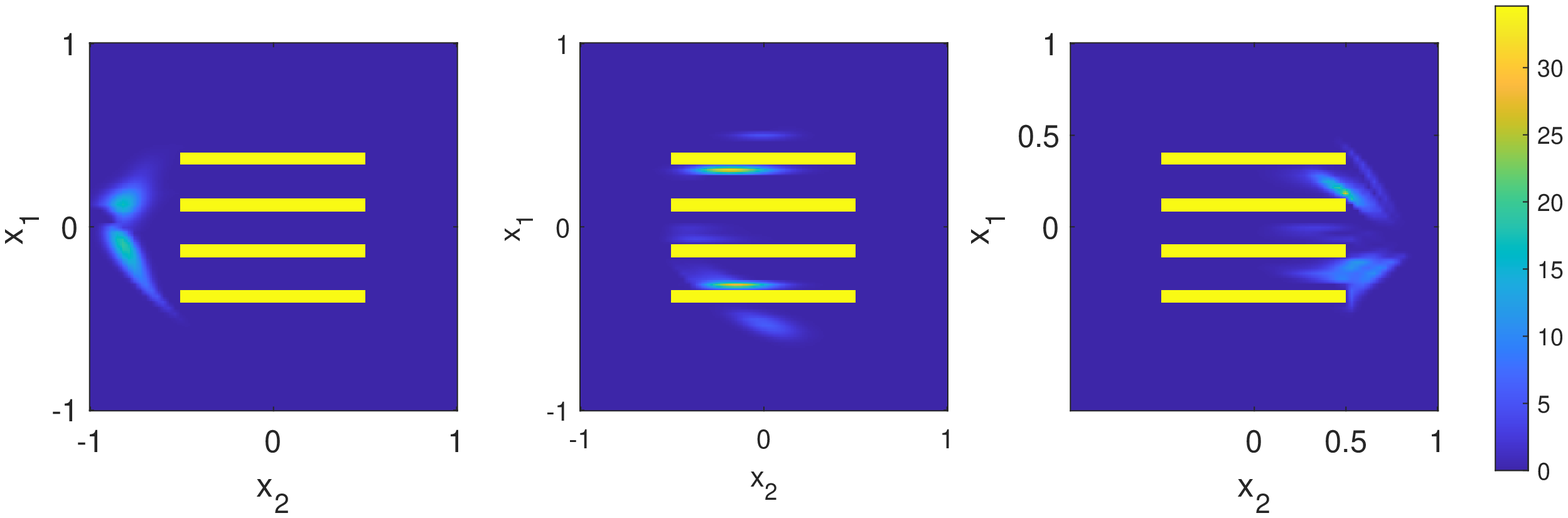}
			\vspace{-0.7cm}
			\caption{$\sigma_1 =0.2, \sigma_2 = 0.5, \mu =5 $}
		\end{subfigure}
		\caption{\textbf{Gaussian repulsion with static obstacles.} MFG solution $\rho(x,t)$ at $t=0.1, 0.5, 0.9$ with different Gaussian parameters $(\sigma_1,\sigma_2,\mu)$, where bright yellow rectangles represents static obstacles.}
		\label{fig:gaussian_static_obstacle}
	\end{figure}

	We consider a MFG model with Gaussian repulsion from Section \ref{subsec:gauss} on the domain $\Omega \times [0,T] = [-1,1]^2 \times [0,1]$. We set the initial-terminal conditions for our MFG system to be
	\begin{equation*}
	\begin{split}
	\rho_0(x_1, x_2) & = \rho_G(0, -0.9,0.2)\\
%	\rho_0(x_1, x_2) = \frac{1}{2 \pi \sigma^2 } \exp \left(-\dfrac{x_1^2 + (x_2+0.9)^2 }{2\sigma^2}\right), \quad \sigma = 0.2\\
	g(x_1, x_2) & = 2\exp\left(-5x_1^2-0.25 
	\left(x_2-0.9\right)^2\right) \left(\left(x_2-0.9 \right)^2 -1\right) +x_1^2\\
	\end{split}
	\end{equation*}	
We fix this $g$ for all examples with Gaussian repulsion. Furthermore, we set
\begin{equation*}
	L(x,v)=\frac{1}{2}\|v\|^2 + 10^3 \left(\max \left(|x_1-0.5|,|x_2-0.5| \right)\right)^8+Q(x),
\end{equation*}
where $Q(x)$ takes on extremely high values on the four rectangular regions in Figure \ref{fig:gaussian_static_obstacle} and 0 elsewhere. Thus, $Q$ models four static rectangular obstacles. Finally, we choose $n=3$ to approximate the kernel.

%	    \begin{figure} [!htbp]
%		\centering
%		\includegraphics[width=1.0\textwidth,trim=0 0 0 0, clip=false]{gaussian_static_obstacle.eps}
%		\vspace{-0.7cm}
%		\caption{MFG solution $\rho(x,0.5)$ for different choice of $\sigma_i$, with bright yellow rectangles represents static obstacles}
%		\label{fig:gaussian_static_obstacle}
%	\end{figure}
 We have computed the MFG solutions for four choices of parameters
\begin{equation*}
(\sigma_1,\sigma_2,\mu) \in \left\{ (0.8,0.8,0.1),(0.2,0.2,5),(0.5,0.2,5),(0.2,0.5,5) \right\}
\end{equation*}
The results are shown in Figure \ref{fig:gaussian_static_obstacle}. As we can see, agents travel through the channels created by $Q(x)$ to avoid high cost. Recall that small $\sigma_i$ yields strong repulsion in $x_i$ direction, which results in different behavior by agents. For instance, in Figure \ref{fig:gaussian_static_obstacle} (C) we impose a strong repulsion in $x_2$ direction and see horizontally elongated density evolution. 
	 
%	 g:2*exp(-5*((x1-0)^2+0.05*(x2-0.9)^2))*((x2-0.9)^2 -1)+ 1*x1^2;
	 %rho0(i,j) = 1/sigma/sqrt(2*pi)*exp(-0.5*(((x2-0.9)^2)/(sigma^2)))* 1/sigma/sqrt(2*pi)*exp(-0.5*(((x1-0.5)^2)/(sigma^2)));

\subsection{Gaussian repulsion with dynamic obstacles}
Next we consider a MFG model with Gaussian repulsion on $\Omega \times [0,T] = [-1,1]^2 \times [0,1]$ with dynamic obstacles. We set
%rho0(i,j) = 1/sigma^2/(2*pi)*exp(-0.5*(((x2+0.9)^2 + (x1-0)^2)/(sigma^2))) + 1/sigma^2/(2*pi)*exp(-0.5*(((x2+0.9)^2 + (x1-0.4)^2)/(sigma^2))) +1/sigma^2/(2*pi)*exp(-0.5*(((x2+0.9)^2 + (x1-0.8)^2)/(sigma^2)))+1/sigma^2/(2*pi)*exp(-0.5*(((x2+0.9)^2 + (x1+0.4)^2)/(sigma^2)))+1/sigma^2/(2*pi)*exp(-0.5*(((x2+0.9)^2 + (x1+0.8)^2)/(sigma^2)));
\begin{equation*}
\rho_0(x) =\frac{1}{5} \sum_{j=1}^5\rho_G(c_j,-0.9,\sigma_G), \quad \sigma_G =0.2, \; c_j = -1.2 + 0.4j,~1\leq j \leq 5
\end{equation*}

		\begin{figure} [!htbp]
		\centering
		\begin{subfigure}{\linewidth}
			\includegraphics[width=1.0\textwidth,trim=0 0 0 0, clip=false]{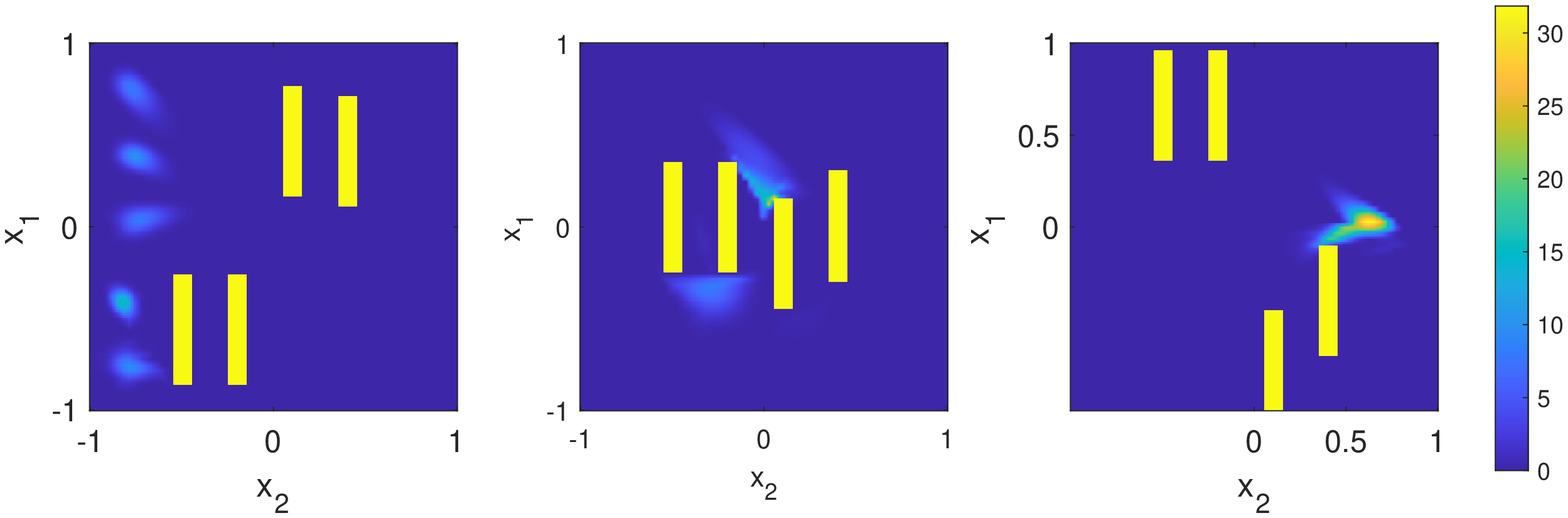}
			\vspace{-0.7cm}
			\caption{$\sigma_1 = 0.8, \sigma_2 = 0.8, \mu = 0.1 $}
		\end{subfigure}
		\begin{subfigure}{\linewidth}
			\includegraphics[width=1.0\textwidth,trim=0 0 0 0, clip=false]{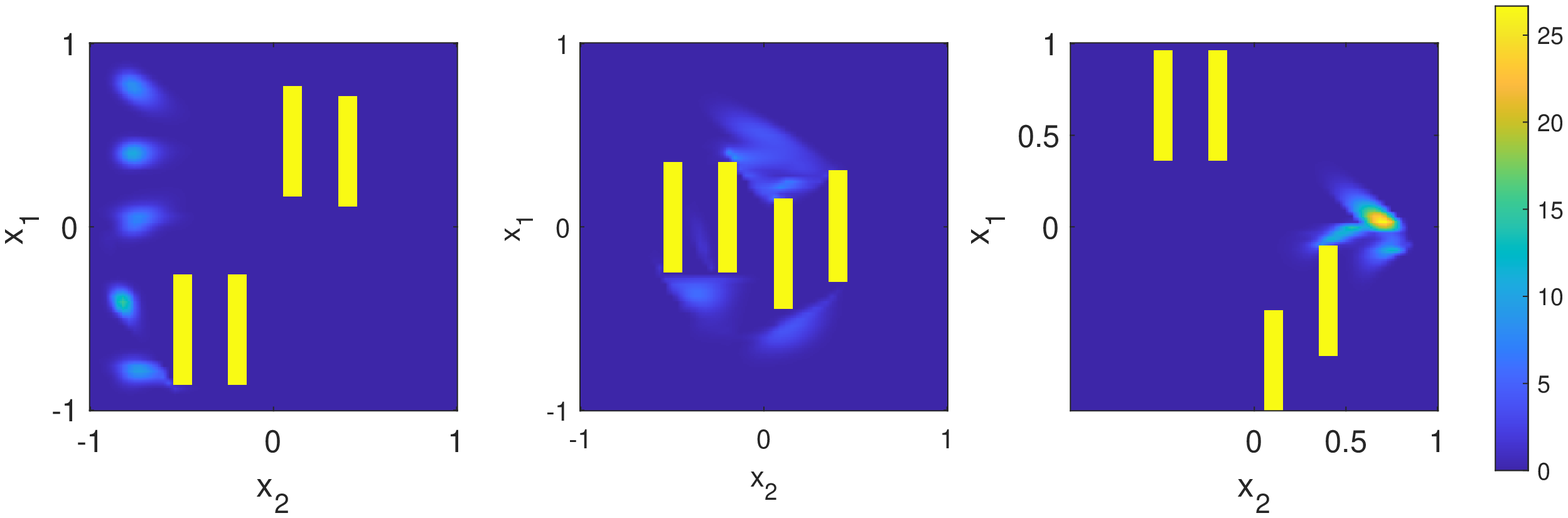}
			\vspace{-0.7cm}
			\caption{$\sigma_1 = 0.2, \sigma_2 = 0.2, \mu = 5$}
		\end{subfigure}
		\begin{subfigure}{\linewidth}
			\includegraphics[width=1.0\textwidth,trim=0 0 0 0, clip=false]{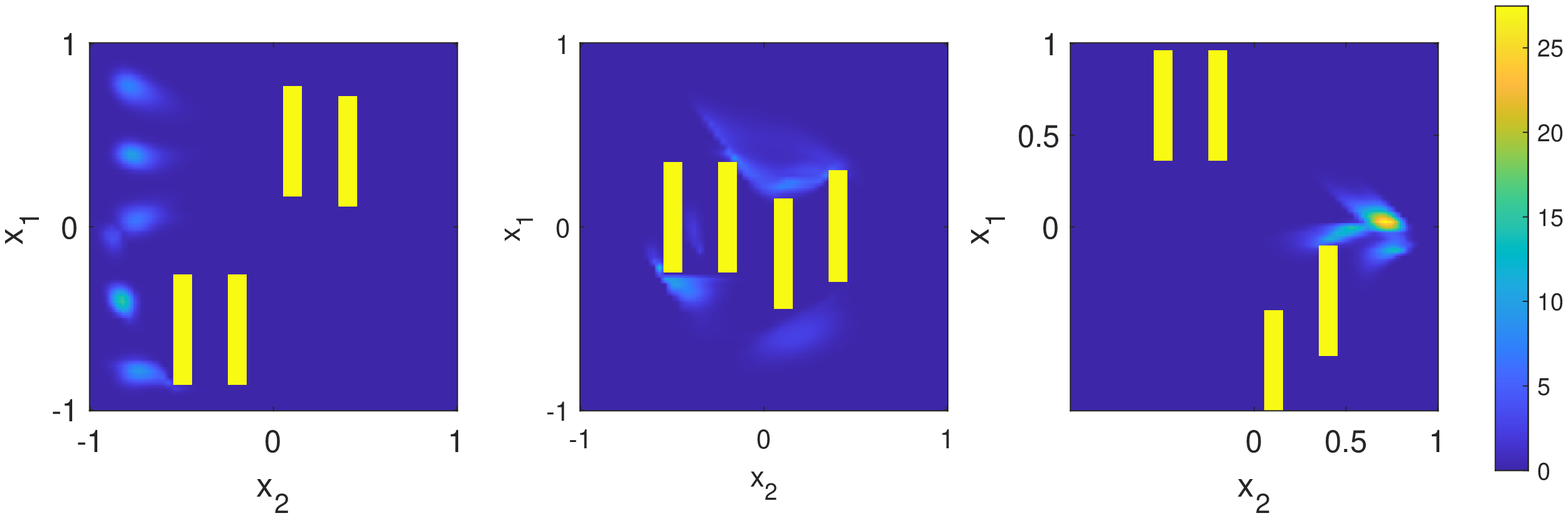}
			\vspace{-0.7cm}
			\caption{$\sigma_1 =0.5, \sigma_2 = 0.2, \mu = 5 $}
		\end{subfigure}
		\begin{subfigure}{\linewidth}
			\includegraphics[width=1.0\textwidth,trim=0 0 0 0, clip=false]{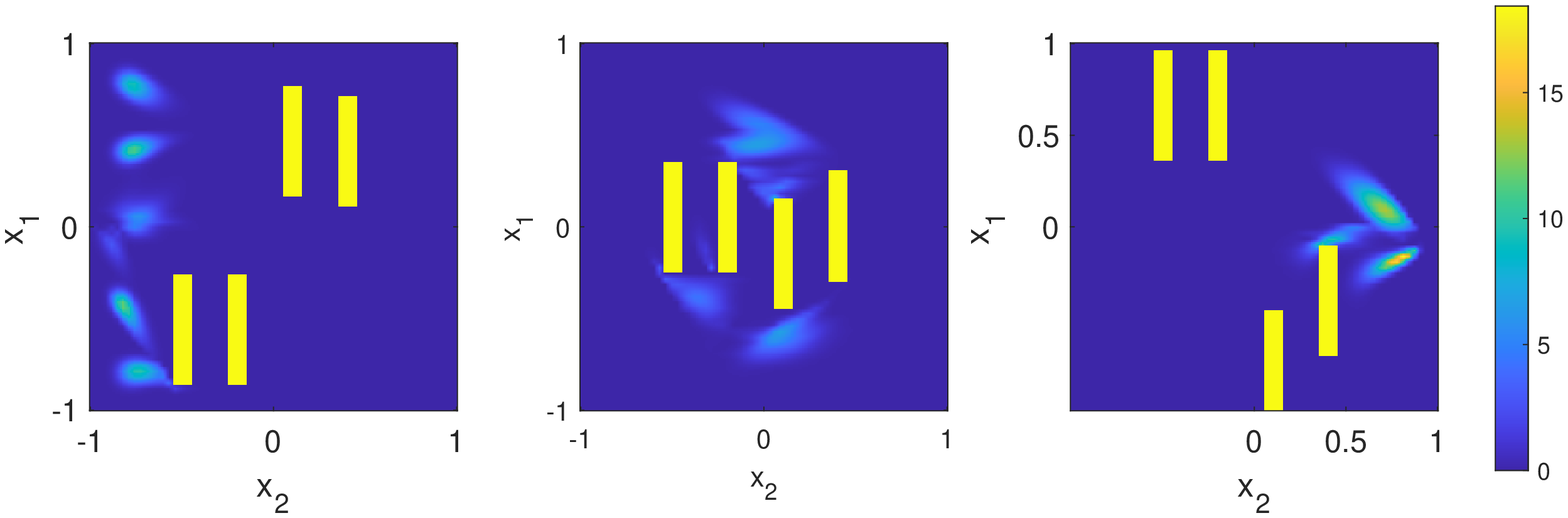}
			\vspace{-0.7cm}
			\caption{$\sigma_1 =0.2, \sigma_2 = 0.5, \mu =5 $}
		\end{subfigure}
		\caption{\textbf{Gaussian repulsion with dynamic obstacles.} MFG solution $\rho(x,t)$ at $t=0.1, 0.5, 0.9$ with different Gaussian parameters $\sigma_1,\sigma_2,\mu$, where bright yellow rectangles represent obstacles moving along $x_1$ direction.} 
		\label{fig:gaussin_moving_obstacle}
	\end{figure} 
%	 \begin{figure} [!htbp]
%		\centering
%		\includegraphics[width=1.0\textwidth,trim=0 0 0 0, clip=false]{gaussian_moving_obstacles.eps}
%		\vspace{-0.7cm}
%		\caption{MFG solution $\rho(x,0.3), \rho(x,0.6)$ for case 1 $\sigma_1=\sigma_2 = 0.8, \mu = 0.1$ (top) and case 2 $\sigma_1=0.5, \sigma_2 = 0.2, \mu = 5$(bottom), with bright yellow rectangles represents obstacles moving along $x_1$ direction.}
%		\label{fig:gaussin_moving_obstacle}
%	\end{figure}
To model dynamic obstacles, we set
\begin{equation*}
L(x,v)=\frac{1}{2}\|v\|^2 + 10^3 \left(\max \left(|x_1-0.5|,|x_2-0.5| \right)\right)^8+Q(x,t)
\end{equation*}
where $Q$ now represents time-dependent rectangular obstacles that move vertically. The rest of the parameters are the same as in the previous section. The results are shown in Figure \ref{fig:gaussin_moving_obstacle}. Again, values of $\sigma_1,\sigma_2$ control how spread is the solution in $x_1,x_2$ directions.

We also note that the computational cost for the static and dynamic models is the same.

\subsection{Interactions in sub-regions}

Next, we consider a MFG model with a Gaussian repulsion on $\Omega \times [0,T] = [-1,1]^2 \times [0,1]$ where agents interact only within domains $\Omega_1=\left\{(x_1,x_2):x_1\leq 0\right\}$ and $\Omega_2=\left\{(x_1,x_2):x_1> 0\right\}$. This means that agents in $\Omega_i$ interact only with those in $\Omega_i$. We set  %rho0(i,j) = 1/sigma^2/(2*pi)*exp(-0.5*(((x2+0.9)^2 + (x1-0)^2)/(sigma^2))) + 1/sigma^2/(2*pi)*exp(-0.5*(((x2+0.9)^2 + (x1-0.4)^2)/(sigma^2))) +1/sigma^2/(2*pi)*exp(-0.5*(((x2+0.9)^2 + (x1-0.8)^2)/(sigma^2)))+1/sigma^2/(2*pi)*exp(-0.5*(((x2+0.9)^2 + (x1+0.4)^2)/(sigma^2)))+1/sigma^2/(2*pi)*exp(-0.5*(((x2+0.9)^2 + (x1+0.8)^2)/(sigma^2)));
\begin{equation*}
\begin{split}
	g(x_1, x_2 ) & = -4 \exp\left(- 5\left(x_1-0.0)^2-2.5(x_2-0.5)^2\right) \right)\\
\rho_0(x) & =\frac{1}{2} \rho_G(0.2,-0.9,\sigma_G) + \frac{1}{2}\rho_G(-0.2,-0.9,\sigma_G) 
\end{split}
\end{equation*}
The rest of the parameters are the same as for previous examples with Gaussian repulsion. We apply the basis modification explained in Section \ref{subsec:subregions} to compute the solution.
%	\begin{figure} [!htbp]
%		\centering
%		\includegraphics[width=1.0\textwidth,trim=0 0 0 0, clip=false]{subregion_gaussian.eps}
%		\vspace{-0.7cm}
%		\caption{MFG solution $\rho(x,0.3), \rho(x,0.8)$ for sub-region and global interactions. $\sigma_1=\sigma_2 = 0.2, \mu 5$.}
%		\label{fig:gaussin_subregion}
%	\end{figure}
	\begin{figure} [!htbp]
	\centering
	\begin{subfigure}{\linewidth}
		\includegraphics[width=1.0\textwidth,trim=0 0 0 0, clip=false]{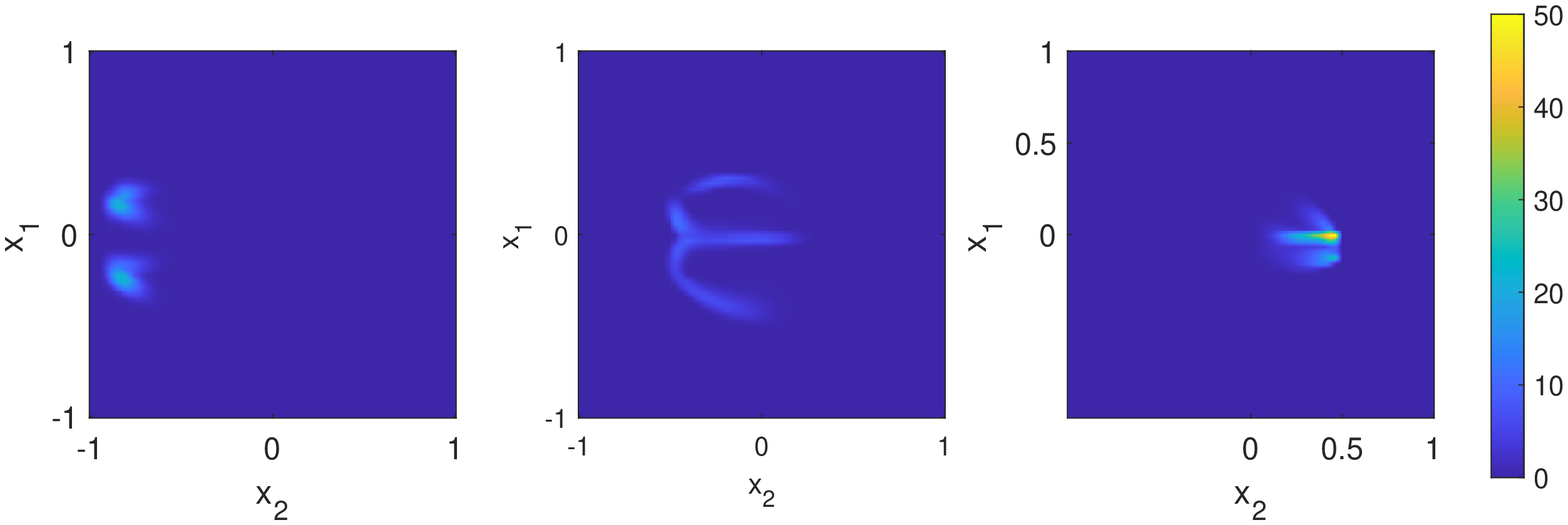}
		\vspace{-0.7cm}
		\caption{sub-region}
	\end{subfigure}
	\begin{subfigure}{\linewidth}
		\includegraphics[width=1.0\textwidth,trim=0 0 0 0, clip=false]{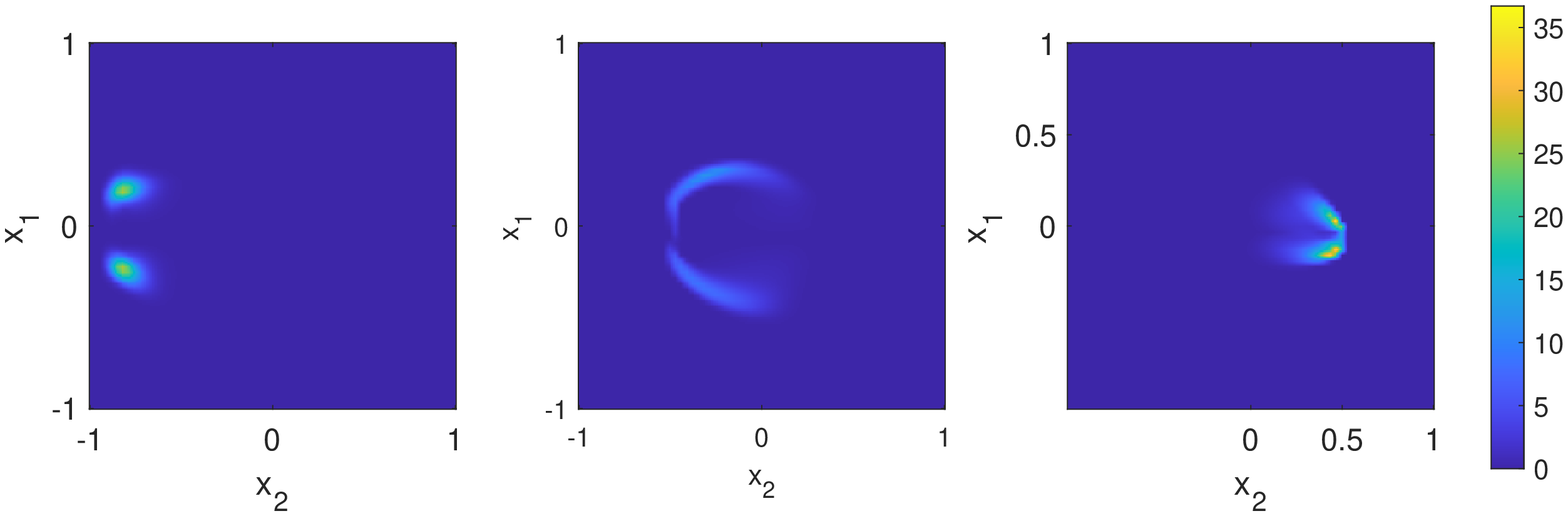}
		\vspace{-0.7cm}
		\caption{global}
	\end{subfigure}
	\caption{\textbf{Sub-region interactions.} MFG solution $\rho(x,0.1), \rho(x,0.5), \rho(x,0.9)$ for sub-region and global interactions with parameters $\sigma_1=\sigma_2 = 0.2,~\mu = 5$}
	\label{fig:gaussin_subregion}
\end{figure} 
The results are shown in Figure \ref{fig:gaussin_subregion} where we have also included the solution with same data but full interaction. In both cases, densities spread before concentrating at the desired location. However, in the sub-region interaction case, Figure \ref{fig:gaussin_subregion} (A), there is a concentration of agents along the common boundary of $\Omega_1,\Omega_2$. The reason is that agents on different sides of this boundary do not interact with each other, so they do not mind congestion.

\subsection{Differential-operator interactions}

Throughout this section, we set $V=\mu (I-\Delta )^{-2},~\mu>0$ and $\Omega=\mathbb{T}^d$. In \cite{achdou10}, the authors solve a stationary MFG system
\begin{equation}\label{eq:mfg_stat}
\begin{cases}
H(x,\nabla \phi) = V[\rho]+\lambda\\
-\nabla \cdot(\rho \nabla_p H(x,\nabla \phi))=0\\
\int \rho=1,~ \rho \geq 0,\quad ~ \lambda \in \mathbb{R}
\end{cases}
\end{equation}
by approximating it with its second-order version
\begin{equation}\label{eq:mfg_achdou}
\begin{cases}
-\sigma \Delta \phi+H(x,\nabla \phi) = V[\rho]+\lambda\\
-\sigma \Delta \rho -\nabla \cdot(\rho \nabla_p H(x,\nabla \phi))=0\\
\int \rho=1,~ \rho \geq 0,\quad ~ \lambda \in \mathbb{R}
\end{cases}	
\end{equation}
for small $\sigma>0$. Here, we recover the results in \cite{achdou10} using our method. Since we consider first-order time-dependent systems instead of second-order stationary ones, we apply a different approximation procedure for \eqref{eq:mfg_stat}. 

To approximate \eqref{eq:mfg_stat} we use the long-time convergence, or the turnpike property, of MFG systems discussed in \cite{carda13}. More precisely, we approximate \eqref{eq:mfg_stat} by
\begin{equation}\label{eq:mfg_T}
\begin{cases}
-\psi_t+H(x,\nabla \psi) = V[\nu]\\
\nu_t-\nabla \cdot(\nu \nabla_p H(x,\nabla \psi))=0\\
\nu(x,0)=\rho_0(x),~ \psi(x,T)=g(x)
\end{cases}
\end{equation}
where $\rho_0 \in L^\infty(\mathbb{T}^d),~g\in C^2(\mathbb{T}^d)$ and $T>0$ is large. To formulate the convergence results, we need to scale the time variable in \eqref{eq:mfg_T} and obtain a problem on a time-interval $[0,1]$. For that, we write $\psi(x,t)=\phi(x,\frac{t}{T}),~\nu(x,t)=\rho(x,\frac{t}{T})$, and \eqref{eq:mfg_T} becomes
\begin{equation}\label{eq:mfg_T=1}
\begin{cases}
-\phi_t+T \cdot H(x,\nabla \phi) = T V[\rho]\\
\rho_t-\nabla \cdot\big(\rho \nabla_p (T \cdot H(x,\nabla \phi))\big)=0\\
\rho(x,0)=\rho_0(x),~ \phi(x,1)=g(x)
\end{cases}
\end{equation}
Furthermore, a triple $(\bar{\phi},\bar{\rho},\bar{\lambda})$ is solution of \eqref{eq:mfg_stat} if $\bar{\phi}$ is a Lipschitz viscosity solution of the HJB in \eqref{eq:mfg_stat}, $\nabla \bar{\phi}$ exists $\bar{\rho}$ a.e., and the continuity equation in \eqref{eq:mfg_stat} is satisfied in the sense of distributions. We summarize the results from \cite{carda13} in the following theorem. We omit assumptions and technicalities and refer to the original paper for details.
\begin{theorem}[\cite{carda13}]\label{thm:carda}
Under suitable assumptions,
\begin{enumerate}
	\item system \eqref{eq:mfg_stat} has at least one solution. Moreover, if $(\bar{\phi}_1,\bar{\rho}_1,\bar{\lambda}_1)$ and $(\bar{\phi}_2,\bar{\rho}_2,\bar{\lambda}_2)$ are solutions, then $\bar{\lambda}_1=\bar{\lambda}_2$, and $V[\bar{\rho}_1]=V[\bar{\rho}_2]$.
	\item for a solution $(\bar{\phi},\bar{\rho},\bar{\lambda})$ of \eqref{eq:mfg_stat} one has that
\end{enumerate}
	\begin{equation*}
	\sup_{t\in [0,1]} \left\| \frac{\phi_T(\cdot,t)}{T}-\bar{\lambda}(1-t)\right\|_{L^\infty(\mathbb{T}^d)} \leq \frac{C}{T^{\frac{1}{2}}},
	\end{equation*}
	and
	\begin{equation*}
	\int_0^1 \|V[\rho_T(\cdot,t)]-V[\bar{\rho}]\|_{L^\infty(\mathbb{T}^d)} dt \leq \frac{C}{T^{\frac{1}{2}}},
	\end{equation*}
	where $(\phi_T,\rho_T)$ is the solution of \eqref{eq:mfg_T=1}.
\end{theorem}
Therefore, to approximate solutions of \eqref{eq:mfg_stat} we need to solve \eqref{eq:mfg_T=1}. We take $H(x,p)=\frac{|p|^2}{2}-Q(x)$ where $Q$ is a smooth periodic function. In this case, one can easily verify that assumptions in \cite{carda13} are fulfilled.

As mentioned in the theorem above, a solution $\bar{\phi}$ in \eqref{eq:mfg_stat} is not necessarily unique even up to constants, whereas $\bar{\lambda}$ and $V[\bar{\rho}]$ are. However, for $V=\mu (I-\Delta)^{-2}$ the uniqueness of $V[\bar{\rho}]$ implies that of $\bar{\rho}$. Furthermore, $\lim\limits_{T\to \infty} \|V[\rho_T(\cdot,t)]-V[\hat{\rho}]\|_{L^\infty(\mathbb{T}^d)}=0$ implies a weak convergence $\rho_T(\cdot,t) \rightharpoonup \bar{\rho}$. Hence Theorem \ref{thm:carda} guarantees that for a large set of times $t\in [0,1]$ the solution $\rho_T(\cdot,t)$ of \eqref{eq:mfg_T=1} converges weakly to a well defined limit $\bar{\rho}$.

\begin{figure} [!htbp]
\centering
\includegraphics[width=1.0\textwidth,trim=0 0 0 0, clip=false]{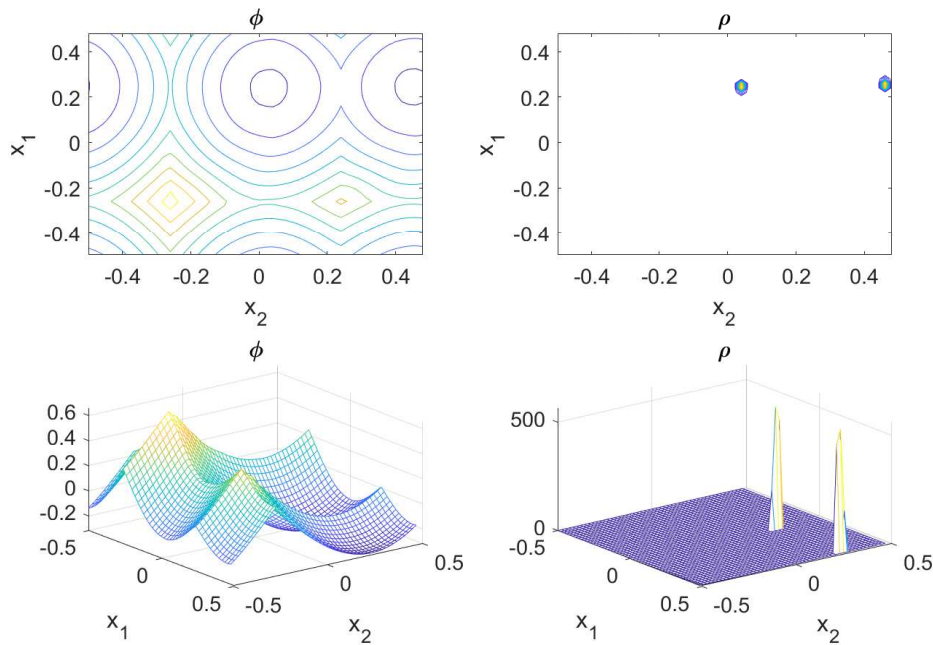}
\vspace{-0.7cm}
\caption{The contours and graphs of $\phi_T(x,t)-\int_{\mathbb{T}^d} \phi_T(y,t)dy$ and $\rho_T(x,t)$ for $T=10$ and $t=0.4$.}
\label{fig:achdou_eg}
\end{figure}

As in \cite{achdou10} we take $d=2,\mu=200$, and
\begin{equation*}
Q(x_1,x_2)=-\sin (2\pi x_2)-\sin (2\pi x_1)-\cos(4\pi x_1)
\end{equation*}
We approximate $V$ as in Section \ref{subsec:diffop} using trigonometric polynomials up to order $n=2$. Additionally, we set $T=10,~\rho_0(x)=1,~g(x)=0$ . The results are shown in Figure \ref{fig:achdou_eg} where we plot $\rho_T(x,t)$. We also plot $\phi_T(x,t)-\int_{\mathbb{T}^d} \phi_T(y,t)dy$ to test whether it approximates a solution $\bar{\phi}$ of \eqref{eq:mfg_stat}. Latter holds for second-order problems but not the first-order ones.

As we can see, we obtain accurate reconstructions of Tests 5, 6 in \cite{achdou10}. Our solutions are slightly less diffused because we consider first-order equations as opposed to second-order ones in \cite{achdou10}. Additionally, we use $H(x,p)=\frac{|p|^2}{2}-Q(x)$ whereas the examples in \cite{achdou10} are computed for $H(x,p)=|p|^{\frac{3}{2}}-Q(x)$. Nevertheless, we believe that qualitative properties and shapes of the solutions do not alter much due to this difference.

%\bibliographystyle{plain}
%\bibliography{main}

\end{document}